\documentclass[11pt]{article}
\usepackage{amsfonts}
\usepackage{amsfonts}
\usepackage{amsfonts}
\usepackage{mathrsfs}
\usepackage{bm}
\usepackage{amssymb,amsmath} % font used for R in Real numbers
 \allowdisplaybreaks

\catcode`!=11
\let\!int\int \def\int{\displaystyle\!int}
\let\!lim\lim \def\lim{\displaystyle\!lim}
\let\!sum\sum \def\sum{\displaystyle\!sum}
\let\!sup\sup \def\sup{\displaystyle\!sup}
\let\!inf\inf \def\inf{\displaystyle\!inf}
\let\!cap\cap \def\cap{\displaystyle\!cap}
\let\!max\max \def\max{\displaystyle\!max}
\let\!min\min \def\min{\displaystyle\!min}
\let\!frac\frac \def\frac{\displaystyle\!frac}
\catcode`!=12

\let\oldsection\section
\renewcommand\section{\setcounter{equation}{0}\oldsection}

\allowdisplaybreaks
\def\pf{\it{Proof.}\rm\quad}

\def\N{\mathbb{N}}

\newtheorem{thm}{Theorem}[section]
\newtheorem{lem}[thm]{Lemma}
\newtheorem{cor}[thm]{Corollary}

\begin{document}
%%%%%%%%%%%%%%%%%%%% title %%%%%%%%%%%%%%%%%%%%%%%%%%%%%%%%%%%%%%%%%%%%%%%%
\title {\bf Some evaluation of parametric Euler sums }
\author{
{Ce Xu\thanks{Corresponding author. Email: xuce1242063253@163.com}}\\[1mm]
\small School of Mathematical Sciences, Xiamen University\\
\small Xiamen
361005, P.R. China}

\date{}
\maketitle \noindent{\bf Abstract } In this paper, by using the method of Contour Integral Representations and the Theorem of Residues and integral representations of series, we discuss the analytic representations of parametric Euler sums that involve harmonic numbers through zeta values and rational function series, either linearly or nonlinearly. Furthermore, we give explicit formulae for several parametric quadratic and cubic sums in terms of zeta values and rational series. Moreover, some interesting new consequences and illustrative examples are considered.
\\[2mm]
\noindent{\bf Keywords} Harmonic number; Euler sum; Riemann zeta function; Hurwitz zeta function.
\\[2mm]
\noindent{\bf AMS Subject Classifications (2010):} 11M06; 11M32.
\section{Introduction}
Throughout this article we will use the following definitions and notations. Let $\N:=\{1,2,3,\ldots\}$ be the set of natural numbers, and $\N_0:=\N\bigcup\{0\},\mathbb{N} \setminus \{1\}:=\{2,3,4,\cdots\}$. In this paper,
harmonic numbers, alternating harmonic numbers and their generalizations are classically defined by $$\ \zeta_n(k):=\sum\limits_{j=1}^n\frac {1}{j^k} ,\ L_{n}(k):=\sum\limits_{j=1}^n\frac{(-1)^{j-1}}{j^k},k\in \N,\eqno(1.1)$$
where $H_n:=\zeta_n(1)=\sum\limits_{j=1}^n\frac {1}{j}$ is classical harmonic number and the empty sum $\zeta_0{(m)}$ is conventionally understood to be zero.
The subject of this paper is Euler sums, which are the infinite sums whose general term is a product of harmonic
numbers and alternating harmonic numbers of index $n$ and a power of $n^{-1}$ or ($(-1)^{n-1}n^{-1}$).
Hence, more generally we can define the Euler sums by the series ([18,19])
$$\sum\limits_{n = 1}^\infty  {\frac{{\prod\limits_{i = 1}^{{m_1}} {\zeta _n^{{q_{_i}}}\left( {{k_i}} \right)\prod\limits_{j = 1}^{{m_2}} {L_n^{{l_j}}\left( {{h_j}} \right)} } }}{{{n^p}}}} ,\ \sum\limits_{n = 1}^\infty  {\frac{{\prod\limits_{i = 1}^{{m_1}} {\zeta _n^{{q_{_i}}}\left( {{k_i}} \right)\prod\limits_{j = 1}^{{m_2}} {L_n^{{l_j}}\left( {{h_j}} \right)} } }}{{{n^p}}}}{{\left( { - 1} \right)}^{n - 1}}, \eqno(1.2)$$
where $m_{1},m_{2},q_{i},k_{i},h_{j},l_{j},p(p\geq 2)$ are positive integers. The quantity $w:=\sum\limits_{i=1}^{m_{1}}(k_{i}q_{i})+\sum\limits_{j=1}^{m_{2}}(h_{j}l_{j})+p$ is called the weight,  the quantity $k: = \sum\limits_{i = 1}^{{m_1}} {{q_i}}  + \sum\limits_{j = 1}^{{m_2}} {{l_j}} $  is called the degree.\\
\indent Apart from the actual evaluation of the series, one of the main questions that one sets out to solve is whether or not a given series can be expressed in terms of a linear rational combination of known constants. When this is the case, we say that the series is reducible to these values.
\\
It has been discovered in the course of the years that many Euler sums admit expressions involving finitely the zeta values, that is to say value of the Riemann zeta function,
$$\zeta(s):=\sum\limits_{n = 1}^\infty {\frac {1}{n^{s}}},\Re(s)>1$$
with positive integer arguments. Note that the alternating Riemann zeta function is defined respectively by
\[\bar \zeta \left( s \right) = \sum\limits_{n = 1}^\infty  {\frac{{{{\left( { - 1} \right)}^{n - 1}}}}{{{n^s}}}}=(1-2^{1-s})\zeta(s) ,\;{\mathop{\Re}\nolimits} \left( s \right) \ge 1.\]\\
\indent For a pair $(p,q)$ of positive integers with $q \geq 2$, the classical linear Euler sum is defined by
$$S_{p,q}:=\sum\limits_{n = 1}^\infty  {\frac{1}{{{n^q}}}} \sum\limits_{k = 1}^n {\frac{1}{{{k^p}}}}. \eqno(1.3)$$
The study of these Euler sums was started by Euler. The earliest results on Euler sums are due
to Euler who elaborated a method to reduce linear sums of small weight to certain rational linear
combinations of products of zeta values. In 1742, Goldbach proposed to Euler the problem of expressing the $S_{p,q}$ in terms of values at positive integers of the Riemann zeta function $\zeta(s)$. Euler showed this problem in the case $p = 1$ and gave a general formula for odd weight $p + q$ in 1775. Moreover, he conjectured that the double linear
sums would be reducible to zeta values when $p + q$ is odd, and even gave what he hoped to obtain the general
formula. In [4], D. Borwein, J.M. Borwein and R. Girgensohn proved conjecture and formula, and in [1], D.H. Bailey, J.M. Borwein and R. Girgensohn conjectured that the double linear sums $S_{p,q}$ when $p + q > 7,p + q$ is even, are not reducible.\\
\indent Let $\pi  = \left( {{\pi _1}, \ldots ,{\pi _k}} \right)$ be a partition of integer $p$ and $p = {\pi _1} +  \cdots  + {\pi _k}$ with ${\pi _1} \le {\pi _2} \le  \cdots  \le {\pi _k}$. The classical nonlinear Euler sum of index $\pi,q$ is defined as follows (see [9])
\[{S_{\pi ,q}} := \sum\limits_{n = 1}^\infty  {\frac{{{\zeta _n}\left( {{\pi _1}} \right){\zeta _n}\left( {{\pi _2}} \right) \cdots {\zeta _n}\left( {{\pi _k}} \right)}}{{{n^q}}}},\]
where the quantity ${\pi _1} +  \cdots  + {\pi _k} + q$ is called the weight, the quantity $k$ is called the degree.
The relationship between the values of the Riemann zeta function and nonlinear Euler sums has been studied by many authors, for example see [1-11,13-19]. Euler sums may be studied through a profusion of methods: combinatorial, analytic and algebraic. Philippe Flajolet and Bruno Salvy informed us about
some ongoing work of theirs ([9]) to evaluate Euler sums in an entirely different way, namely
using contour integration and the residue theorem. In this way they manage to prove, for
example, that the cubic sums
\[{S_{{1^3},q}} = \sum\limits_{n = 1}^\infty  {\frac{{H_n^3}}{{{n^q}}}} \;\left( {q = 2,3,4,6} \right)\]
can be evaluated in terms of Riemann zeta values. Furthermore, they proved the quadratic sums
\[{S_{{p_1}{p_2},q}} = \sum\limits_{n = 1}^\infty  {\frac{{{\zeta _n}\left( {{p_1}} \right){\zeta _n}\left( {{p_2}} \right)}}{{{n^q}}}} \]
can be expresses as a rational linear combination of products of linear sums and zeta values whenever $p_1+p_2+q$ is even, and $p_1>1, p_2>1, q>1$. In [19], we showed that all quadratic Euler sums of the form
\[S_{1m,p}=\sum\limits_{n = 1}^\infty  {\frac{{{H_n}{\zeta_n{(m)} }}}{{{n^p}}}}\ \ \left( 3\leq {m + p \le 8} \right)\]
can be reduced to polynomials in zeta values and linear sums.\\
\indent So far, surprisingly little work has been done on parametric
Euler sums. Similarly to the definition of (1.3), the parametric linear Euler sum is defined by the series
\[{S_{p,q}}\left( {{a_1}, \cdots ,{a_p}} \right) := \sum\limits_{n = 1}^\infty  {\frac{{{\zeta _n}\left( p \right)}}{{R_n\left( a_1,a_2,\ldots,a_q \right) }}},\eqno(1.4)\]
\[R_n\left( a_1,a_2,\ldots,a_q \right) := {n^q} + {a_1}{n^{q - 1}} +  \cdots  + {a_q},\]
where ${R_n\left( a_1,a_2,\ldots,a_q \right)}\neq 0$ is a rational function and $ deg({R_n\left( a_1,a_2,\ldots,a_q \right)}=q\geq 2$ .\\
\indent Similarly to the definition of (1.2), the generalized parametric Euler sums is defined as
$$\sum\limits_{n = 1}^\infty  {\frac{{\prod\limits_{i = 1}^{{m_1}} {\zeta _n^{{q_{_i}}}\left( {{k_i}} \right)\prod\limits_{j = 1}^{{m_2}} {L_n^{{l_j}}\left( {{h_j}} \right)} } }}{{{R_n\left( a_1,a_2,\ldots,a_q \right)}}}} ,\ \sum\limits_{n = 1}^\infty  {\frac{{\prod\limits_{i = 1}^{{m_1}} {\zeta _n^{{q_{_i}}}\left( {{k_i}} \right)\prod\limits_{j = 1}^{{m_2}} {L_n^{{l_j}}\left( {{h_j}} \right)} } }}{{{R_n\left( a_1,a_2,\ldots,a_q \right)}}}}{{\left( { - 1} \right)}^{n - 1}}.
\eqno(1.5)$$
In the paper, we will consider the following type of parametric linear sums involving harmonic numbers
\[\sum\limits_{n = 1}^\infty  {\frac{{{H_n}}}{{{{\left( {n + a} \right)}^{s + 1}}}}} ,\;\sum\limits_{n = 1}^\infty  {\frac{{{\zeta _n}\left( {2m - 1} \right)}}{{{n^{2s}}\left( {{n^2} - {a^2}} \right)}}} ,\;\sum\limits_{n = 1}^\infty  {\frac{{{\zeta _n}\left( {2m} \right)}}{{n\left( {{n^2} - {a^2}} \right)}}}.\eqno(1.6)\]
and parametric quadratic, cubic Euler sums of the form
\[\sum\limits_{n = 1}^\infty  {\frac{{H_n^2}}{{{{\left( {n + b} \right)}^{s + 1}}}}} ,\;\sum\limits_{n = 1}^\infty  {\frac{{H_n^3}}{{{{\left( {n + b} \right)}^{s + 1}}}}}. \eqno(1.7)\]
where $m,s$ are positive integers and $a \ne  \pm 1, \pm 2, \cdots ,\;b \ne  - 1, - 2, \cdots .$\\
We prove that the sums of (1.6) can be expressed as a rational linear combination of several given rational series and
the parametric quadratic Euler sums of (1.7) are reducible to parametric linear sums.
For example, we prove the identity
\begin{align*}
\sum\limits_{n = 1}^\infty  {\frac{{{H_n}}}{{{{\left( {n + a} \right)}^s}}}}
 &=\frac{s}{2}\zeta \left( {s + 1,a + 1} \right) - \frac{1}{2}\sum\limits_{j = 1}^{s - 2} {\zeta \left( {s - j,a + 1} \right)} \zeta \left( {j + 1,a + 1} \right)
\nonumber \\
           &\quad + a\zeta \left( {s,a + 1} \right)\sum\limits_{n = 1}^\infty  {\frac{1}{{n\left( {n + a} \right)}}}  + \sum\limits_{n = 1}^\infty  {\frac{1}{{n{{\left( {n + a} \right)}^s}}}},\tag{1.8}
\end{align*}
where $s \in \mathbb{N}\setminus\{1\}$.
 $\zeta \left( {s,a + 1} \right)$ and $\bar \zeta \left( {s,a + 1} \right)$ stand for the Hurwitz zeta function and alternating Hurwitz zeta function defined by
\[\zeta \left( {s,a + 1} \right) = \sum\limits_{n = 1}^\infty  {\frac{1}{{{{\left( {n + a} \right)}^s}}}} ,\bar \zeta \left( {s,a + 1} \right) = \sum\limits_{n = 1}^\infty  {\frac{{{{\left( { - 1} \right)}^{n - 1}}}}{{{{\left( {n + a} \right)}^s}}}} ,\;{\mathop{\Re}\nolimits}(s) > 1,\;a \ne  - 1, - 2, \cdots. \tag{1.9}\]
Similarly, the parametric harmonic number (also called partial sums of Hurwitz zeta function) ${{\zeta _n}\left( {p,a} \right)}$ for $p\geq1$ is defined as
\[{\zeta _n}\left( {p,a+1} \right) = \sum\limits_{k = 1}^n {\frac{1}{{{{\left( {k + a} \right)}^p}}}},a\neq -1,-2,\ldots.\]

\section{Parametric linear Euler sums}
In this section we consider the following type of parametric linear Euler sums involving harmonic numbers by the method of contour integration
\[\sum\limits_{n = 1}^\infty  {\frac{{{\zeta _n}\left( {2m - 1} \right)}}{{{n^{2s}}\left( {{n^2} - {a^2}} \right)}}} ,\;\sum\limits_{n = 1}^\infty  {\frac{{{\zeta _n}\left( {2m} \right)}}{{n\left( {{n^2} - {a^2}} \right)}}}, s\in\N_0.\]
Contour integration is a classical technique for evaluating infinite sums by reducing them to a finite number of residue computations. This summation mechanism is formalized by a lemma that goes back to Cauchy and is nicely developed throughout [9]. Next, we give two lemmas. The following lemma will be useful in the
development of the main theorems.
\begin{lem}([9])
Let $\xi \left( s \right)$ be a kernel function and let $r(s)$ be a rational function which is $O(s^{-2})$ at infinity. Then
\[\sum\limits_{\alpha  \in O} {{\mathop{\rm Re}\nolimits} s{{\left( {r\left( s \right)\xi \left( s \right)} \right)}_{s = \alpha }}}  + \sum\limits_{\beta  \in S} {{\mathop{\rm Re}\nolimits} s{{\left( {r\left( s \right)\xi \left( s \right)} \right)}_{s = \beta }}}  = 0.\tag{2.1}\]
where $S$ is the set of poles of $r(s)$ and $O$ is the set of poles of $\xi \left( s \right)$ that are not poles $r(s)$ . Here ${\mathop{\rm Re}\nolimits} s{\left( {r\left( s \right)} \right)_{s = \alpha }} $ denotes the residue of $r(s)$ at $s= \alpha$. The kernel function $\xi \left( s \right)$ is meromorphic in the whole complex plane and satisfies $\xi \left( s \right)=o(s)$ over an infinite collection of circles $\left| z \right| = {\rho _k}$ with ${\rho _k} \to \infty . $
\end{lem}
\begin{lem} For integer $n\in \N$, then the following relations holds
\begin{align*}
&\sum\limits_{j = 1}^n {\frac{{s\left( {j,p} \right)}}{{j!}}}  = \frac{{s\left( {n + 1,p + 1} \right)}}{{n!}},\;p\in \N, \tag{2.2}\\
&\sum\limits_{j = 1}^{n - 1} {\frac{{{H_j}s\left( {n - j,p - 1} \right)}}{{\left( {n - j} \right)!}}}  = p\frac{{s\left( {n + 1,p + 1} \right)}}{{n!}},\;p\in \mathbb{N} \setminus \{1\}.\tag{2.3}
\end{align*}
where ${s\left( {n,k} \right)}$ is called (unsigned) Stirling number of the first kind ([12]) defined by
\[n!\left( {1 + x} \right)\left( {1 + \frac{x}{2}} \right) \cdots \left( {1 + \frac{x}{n}} \right) = \sum\limits_{k = 0}^n {s\left( {n + 1,k + 1} \right){x^k}} ,\]
when $n<k$, $s(n,k)=0$; when $n,k>1$, $s(n,0)=s(0,k)=0$; when $n=k=0$, $s(0,0)=1$.
\end{lem}
\pf From [12], we have the generating function of (unsigned) Stirling number of the first kind:
\[{\ln ^{p+1}}\left( {1 - x} \right) = {\left( { - 1} \right)^{p+1}}(p+1)!\sum\limits_{n = p+1}^\infty  {s\left( {n,p+1} \right)} \frac{{{x^n}}}{{n!}},\ p\in \N_0, -1\leq x<1.\tag{2.4}\]
Differentiating this equality, we obtain
\[\frac{{{{\ln }^{p }}\left( {1 - x} \right)}}{{1 - x}} = {\left( { - 1} \right)^{p}}p!\sum\limits_{n = p }^\infty  {s\left( {n + 1,p+1} \right)} \frac{{{x^n}}}{{n!}}\:\:,p \in \N.\tag{2.5}\]
On the other hand, using the cauchy product of power series, we have
\[\frac{{{{\ln }^p}\left( {1 - x} \right)}}{{1 - x}} = {\left( { - 1} \right)^p}p!\sum\limits_{n = p}^\infty  {\frac{{s\left( {n,p} \right)}}{{n!}}{x^n}} \sum\limits_{n = 1}^\infty  {{x^{n - 1}}}  = {\left( { - 1} \right)^p}p!\sum\limits_{n = 1}^\infty  {\left( {\sum\limits_{j = 1}^n {\frac{{s\left( {j,p} \right)}}{{j!}}} } \right){x^n}}.\tag{2.6} \]
Thus, comparing the coefficients of $x^n$ in (2.5) and (2.6), we can deduce (2.2).
Similarly, using the cauchy product of power series again, we get
\begin{align*}
\frac{{{{\ln }^p}\left( {1 - x} \right)}}{{1 - x}} =& \frac{{\ln \left( {1 - x} \right)}}{{1 - x}}{\ln ^{p - 1}}\left( {1 - x} \right)\\
 =&  - \sum\limits_{n = 1}^\infty  {{H_n}{x^n}} {\left( { - 1} \right)^{p - 1}}\left( {p - 1} \right)!\sum\limits_{n = p - 1}^\infty  {\frac{{s\left( {n,p - 1} \right)}}{{n!}}{x^n}} \\
  =& {\left( { - 1} \right)^p}\left( {p - 1} \right)!\sum\limits_{n = 1}^\infty  {\left( {\sum\limits_{j = 1}^n {\frac{{{H_j}s\left( {n - j + 1,p - 1} \right)}}{{\left( {n - j + 1} \right)!}}} } \right){x^{n + 1}}} .\tag{2.7}
\end{align*}
Combining (2.5) and (2.7), we can obtain (2.3). The proof of Lemma 2.2 is thus completed.\hfill$\square$\\
Moreover, by the definition ${s\left( {n,k} \right)}$, we can rewrite it as
\begin{align*}
\sum\limits_{k = 0}^n {s\left( {n + 1,k + 1} \right){x^k}}  &= n!\exp \left\{ {\sum\limits_{j = 1}^n {\ln \left( {1 + \frac{x}{j}} \right)} } \right\}\\
& = n!\exp \left\{ {\sum\limits_{j = 1}^n {\sum\limits_{k = 1}^\infty  {{{\left( { - 1} \right)}^{k - 1}}\frac{{{x^k}}}{{k{j^k}}}} } } \right\}\\
& = n!\exp \left\{ {\sum\limits_{k = 1}^\infty  {{{\left( { - 1} \right)}^{k - 1}}\frac{{{\zeta _n}\left( k \right){x^k}}}{k}} } \right\}.
\end{align*}
Therefore, we know that ${s\left( {n,k} \right)}$ is a rational linear combination of products of harmonic numbers.
The following identities is easily derived
\begin{align*}
& s\left( {n,1} \right) = \left( {n - 1} \right)!,S\left( {n,2} \right) = \left( {n - 1} \right)!{H_{n - 1}},S\left( {n,3} \right) = \frac{{\left( {n - 1} \right)!}}{2}\left[ {H_{n - 1}^2 - {\zeta _{n - 1}}\left( 2 \right)} \right],\\
&s\left( {n,4} \right) = \frac{{\left( {n - 1} \right)!}}{6}\left[ {H_{n - 1}^3 - 3{H_{n - 1}}{\zeta _{n - 1}}\left( 2 \right) + 2{\zeta _{n - 1}}\left( 3 \right)} \right], \\
&s\left( {n,5} \right) = \frac{{\left( {n - 1} \right)!}}{{24}}\left[ {H_{n - 1}^4 - 6{\zeta _{n - 1}}\left( 4 \right) - 6H_{n - 1}^2{\zeta _{n - 1}}\left( 2 \right) + 3\zeta _{n - 1}^2\left( 2 \right) + 8H_{n - 1}^{}{\zeta _{n - 1}}\left( 3 \right)} \right].
\end{align*}
Therefore, putting $p=2$ in (2.2) and (2.3), we can give the identities
\[\sum\limits_{j = 1}^{n - 1} {\frac{{{H_j}}}{{n - j}}}  = H_n^2 - {\zeta _n}\left( 2 \right),\;\sum\limits_{k = 1}^n {\frac{{{H_k}}}{k}}  = \frac{{H_n^2 + {\zeta _n}\left( 2 \right)}}{2}.\tag{2.8}\]
\indent We make here an essential use of kernels involving the $\psi$ function. The $\psi$ function is the logarithmic derivative of the Gamma function,
\[\psi \left( s \right) = \frac{d}{{ds}}\ln \Gamma \left( s \right) =  - \gamma  - \frac{1}{s} + \sum\limits_{n = 1}^\infty  {\left( {\frac{1}{n} - \frac{1}{{n + s}}} \right)}\tag{2.9} \]
and it satisfies the complement formula
\[\psi \left( s \right) - \psi \left( { - s} \right) =  - \frac{1}{s} - \pi \cot \left( {\pi s} \right),\tag{2.10}\]
as well as an expansion at $s=0$ that involves the zeta values:
\[\psi \left( { - s} \right) + \gamma  = \frac{1}{s} - \sum\limits_{k = 2}^\infty  {\zeta \left( k \right){s^{k - 1}}} .\tag{2.11}\]
Using differentiation $n$ times, we obtain
\[\frac{{{\psi ^{\left( n \right)}}\left( { - s} \right)}}{{n!}} = \frac{1}{{{s^{n + 1}}}} - {\left( { - 1} \right)^n}\sum\limits_{k = 1}^\infty  {\left( {\begin{array}{*{20}{c}}
   {n + k - 1}  \\
   {k - 1}  \\
\end{array}} \right)\zeta \left( {k + n} \right){s^{k - 1}}}, s \to 0
.\tag{2.12}\]
In their paper, ¡°Euler Sums and Contour Integral Representations¡±, Philippe Flajolet and Bruno Salvy gave the following formulae \\
\hrule
\hrule
\begin{align*}
&\pi \cot \left( {\pi s} \right)\mathop  = \limits^{s \to n} \frac{1}{{s - n}} - 2\sum\limits_{k = 1}^\infty  {\zeta \left( {2k} \right){{\left( {s - n} \right)}^{2k - 1}}} , \\
&\psi \left( { - s} \right) + \gamma \mathop  = \limits^{s \to n} \frac{1}{{s - n}} + {H_n} + \sum\limits_{k = 1}^\infty  {\left( {{{\left( { - 1} \right)}^k}{\zeta _n}\left( {k + 1} \right) - \zeta \left( {k + 1} \right)} \right){{\left( {s - n} \right)}^k}} ,\;\;n \ge 0 \\
&\psi \left( { - s} \right) + \gamma \mathop  = \limits^{s \to  - n} {H_{n - 1}} + \sum\limits_{k = 1}^\infty  {\left( {{\zeta _{n - 1}}\left( {k + 1} \right) - \zeta \left( {k + 1} \right)} \right){{\left( {s + n} \right)}^k}} ,\;n > 0\\
&\frac{{{\psi ^{\left( {p - 1} \right)}}\left( { - s} \right)}}{{\left( {p - 1} \right)!}}\mathop  = \limits^{s \to n} \frac{1}{{{{\left( {s - n} \right)}^p}}}\left( {1 + {{\left( { - 1} \right)}^p}\sum\limits_{i \ge p} {\left( {\begin{array}{*{20}{c}}
   {i - 1}  \\
   {p - 1}  \\
\end{array}} \right)\left( {\zeta \left( i \right) + {{\left( { - 1} \right)}^i}{\zeta _n}\left( i \right)} \right){{\left( {s - n} \right)}^i}} } \right),\;n \ge 0,\;p > 1 \\
&\frac{{{\psi ^{\left( {p - 1} \right)}}\left( { - s} \right)}}{{\left( {p - 1} \right)!}}\mathop  = \limits^{s \to  - n} {\left( { - 1} \right)^p}\sum\limits_{i \ge 0} {\left( {\begin{array}{*{20}{c}}
   {p - 1 + i}  \\
   {p - 1}  \\
\end{array}} \right)\left( {\zeta \left( {p + i} \right) - {\zeta _{n - 1}}\left( {p + i} \right)} \right){{\left( {s + n} \right)}^i}} ,\;n > 0,\;p > 1.
\end{align*}
\hrule
\begin{center} Table 1. Local expansions of basic kernels
\end{center}

\indent Nielsen [14], elaborating on Euler's work, proved by a method based on partial fraction expansions that every linear sum $S_{p,q}$ whose weight $p+q$ is odd is expressible as a polynomial in zeta values. We give explicit formula for several classes of parametric Euler sums in terms of Riemann zeta values and rational function series.
Next we evaluate the sums in $(1.6)$.
\begin{thm} Let $m$ be positive integers with $a$ is a real and $a \ne 0, - 1, - 2, \cdots$. Then the following parametric linear sums are reducible to zeta values and rational function series,
\begin{align*}
\sum\limits_{n = 1}^\infty  {\frac{{{\zeta _n}\left( {2m} \right)}}{{n\left( {{n^2} - {a^2}} \right)}}}
 &= \frac{1}{2}\sum\limits_{n = 1}^\infty  {\frac{1}{{{n^{2m + 1}}\left( {{n^2} - {a^2}} \right)}}}  - \frac{{\pi \cot \left( {\pi a} \right)}}{{4{a^2}}}\sum\limits_{n = 1}^\infty  {\left\{ {\frac{1}{{{{\left( {n - a} \right)}^{2m}}}} - \frac{1}{{{{\left( {n + a} \right)}^{2m }}}}} \right\}}
\nonumber \\
& \ \ + \frac{1}{{2{a^2}}}\sum\limits_{k = 1}^m {\zeta \left( {2k} \right)} \sum\limits_{n = 1}^\infty  {\left\{ {\frac{1}{{{{\left( {n + a} \right)}^{2m - 2k + 1}}}} + \frac{1}{{{{\left( {n - a} \right)}^{2m - 2k + 1}}}} - \frac{2}{{{n^{2m - 2k + 1}}}}} \right\}}
\nonumber \\
& \ \ +\frac{m}{{{a^2}}}\zeta \left( {2m + 1} \right) - \frac{1}{{4{a^2}}}\sum\limits_{n = 1}^\infty  {\left\{ {\frac{1}{{{{\left( {n + a} \right)}^{2m + 1}}}} + \frac{1}{{{{\left( {n - a} \right)}^{2m + 1}}}} - \frac{2}{{{n^{2m + 1}}}}} \right\}}
.\tag{2.13}
\end{align*}
\end{thm}
\pf The theorem results from applying the kernel function
\[\xi \left( z \right) = \frac{{\pi \cot \left( {\pi z} \right){\psi ^{\left( {2m - 1} \right)}}\left( { - z} \right)}}{{\left( {2m - 1} \right)!}}\]
to the base function $r\left( z \right) = z^{-1}\left( {{z^2} - {a^2}} \right)^{-1}$. The only singularities are poles at the integers and $\pm a$. At a negative integer $-n$ the pole is simple and the residue is
\[\frac{{{\zeta _n}\left( {2m} \right) - \zeta \left( {2m} \right)}}{{n\left( {{n^2} - {a^2}} \right)}} - \frac{1}{{{n^{2m + 1}}\left( {{n^2} - {a^2}} \right)}}.\]
At a positive integer $n$, the pole has order $2m+1$ and the residue is
\begin{align*}
\frac{{{\zeta _n}\left( {2m} \right) + \zeta \left( {2m} \right)}}{{n\left( {{n^2} - {a^2}} \right)}}&+ \frac{1}{{2{a^2}}}\left\{ {\frac{1}{{{{\left( {n + a} \right)}^{2m + 1}}}} + \frac{1}{{{{\left( {n - a} \right)}^{2m + 1}}}} - \frac{2}{{{n^{2m + 1}}}}} \right\}\\
 &- \frac{1}{{{a^2}}}\sum\limits_{k = 1}^m {\zeta \left( {2k} \right)} \left\{ {\frac{1}{{{{\left( {n + a} \right)}^{2m - 2k + 1}}}} + \frac{1}{{{{\left( {n - a} \right)}^{2m - 2k + 1}}}} - \frac{2}{{{n^{2m - 2k + 1}}}}} \right\}.\\
\end{align*}
The residue of the pole at $\pm a$ is
\[\frac{{\pi \cot \left( {\pi a} \right)}}{{2{a^2}}}\sum\limits_{n = 1}^\infty  {\left\{ {\frac{1}{{{{\left( {n - a} \right)}^{2m}}}} - \frac{1}{{{{\left( {n + a} \right)}^{2m}}}}} \right\}}. \]
Finally the residue of the pole of order $2m+2$ at 0 is found to be
\[ - \frac{{2m}}{{{a^2}}}\zeta \left( {2m + 1} \right).\]
Summing these four contributions yields the statement of the theorem.\hfill$\square$
\begin{thm} Let $m,s\geqslant 0$ be integers with $a$ is a real and $a \ne  0, - 1, - 2, \cdots$. The following parametric linear sums are reducible to zeta values and rational function series,
\begin{align*}
\sum\limits_{n = 1}^\infty  {\frac{{{\zeta _n}\left( {2m + 1} \right)}}{{{n^{2s}}\left( {{n^2} - {a^2}} \right)}}}
 &=\frac{1}{2}\sum\limits_{n = 1}^\infty  {\frac{1}{{{n^{2s + 2m + 1}}\left( {{n^2} - {a^2}} \right)}}}
\nonumber \\
&\quad+ \sum\limits_{n = 1}^s {\sum\limits_{k = 1}^n {\left( {\begin{array}{*{20}{c}}
   {2m + 2k - 1}  \\
   {2k - 1}  \\
\end{array}} \right)\frac{{\zeta \left( {2m + 2k - 1} \right)\zeta \left( {2n - 2k + 2} \right)}}{{{a^{2s - 2n + 2}}}}} }
\nonumber \\
& \quad + \zeta \left( {2m + 1} \right)\left( {\sum\limits_{n = 1}^\infty  {\frac{1}{{{n^2} - {a^2}}}} \left( {\frac{1}{{{n^{2s}}}} - \frac{1}{{{a^{2s}}}}} \right)} \right) \\
&\quad- \frac{1}{2}\sum\limits_{k = 2}^{s + 1} {\left( {\begin{array}{*{20}{c}}
   {2m + 2k - 2}  \\
   {2k - 2}  \\
\end{array}} \right)} \frac{{\zeta \left( {2m + 2k - 1} \right)}}{{{a^{2s - 2k + 4}}}}
\nonumber \\
& \quad+ \frac{{\pi \cot \left( {\pi a} \right)}}{{4{a^{2s + 1}}}}\sum\limits_{n = 1}^\infty  {\left\{ {\frac{1}{{{{\left( {n - a} \right)}^{2m + 1}}}} + \frac{1}{{{{\left( {n + a} \right)}^{2m + 1}}}} - \frac{2}{{{n^{2m + 1}}}}} \right\}}
\nonumber \\
 &\quad- \frac{1}{2}\sum\limits_{j = 1}^s {\frac{{\zeta \left( {2m + 2j + 1} \right)}}{{{a^{2s + 2 - 2j}}}}} \left( {\begin{array}{*{20}{c}}
   {2m + 2j}  \\
   {2j - 1}  \\
\end{array}} \right)
\nonumber \\
& \quad + \sum\limits_{k = 1}^m {\sum\limits_{j = 1}^s {\frac{{\zeta \left( {2k} \right)}}{{{a^{2s + 2 - 2j}}}}} } \left( {\begin{array}{*{20}{c}}
   {2m - 2k + 2j}  \\
   {2j - 1}  \\
\end{array}} \right)\zeta \left( {2m + 2j - 2k + 1} \right)
\nonumber \\
& \quad
 - \frac{1}{{2{a^{2s + 1}}}}\sum\limits_{k = 0}^m {\zeta \left( {2k} \right)} \sum\limits_{n = 1}^\infty  {\left\{ {\frac{1}{{{{\left( {n - a} \right)}^{2m - 2k + 2}}}} - \frac{1}{{{{\left( {n + a} \right)}^{2m - 2k + 2}}}}} \right\}}
,\tag{2.14}
\end{align*}
where the value $\zeta(0)=-\frac{1}{2}$ should be used and $\zeta(1)$ should be replaced by 0 whenever it occurs.
\end{thm}
\pf Similarly to the proof of Theorem 2.3. The theorem results from applying the kernel function
\[\xi \left( z \right) = \frac{{\pi \cot \left( {\pi z} \right){\psi ^{\left( {2m} \right)}}\left( { - z} \right)}}{{\left( {2m} \right)!}}\]
to the base function $r\left( z \right) = {z^{ - 2s}}{\left( {{z^2} - {a^2}} \right)^{ - 1}}$. Note that $r(z)$ can be rewritten as
\[r\left( z \right) = \frac{1}{{{z^{2s}}\left( {{z^2} - {a^2}} \right)}} = \frac{1}{{2{a^{2s + 1}}}}\left( {\frac{1}{{z - a}} - \frac{1}{{z + a}}} \right) - \sum\limits_{j = 1}^s {\frac{1}{{{a^{2s + 2 - 2j}}}} \cdot \frac{1}{{{z^{2j}}}}} .\tag{2.15}\]
From (2.12), (2.15) and Table 1, we can find that, at a positive integer $n$, the pole has order $2m+2$ and the residue is
\begin{align*}
&\frac{{{\zeta _n}\left( {2m + 1} \right) - \zeta \left( {2m + 1} \right)}}{{{n^{2s}}\left( {{n^2} - {a^2}} \right)}}
+ \frac{1}{{{a^{2s + 1}}}}\sum\limits_{k = 1}^m {\zeta \left( {2k} \right)} \left\{ {\frac{1}{{{{\left( {n - a} \right)}^{2m - 2k + 2}}}} - \frac{1}{{{{\left( {n + a} \right)}^{2m - 2k + 2}}}}} \right\} \\
&- 2\sum\limits_{k = 1}^m {\sum\limits_{j = 1}^s {\frac{{\zeta \left( {2k} \right)}}{{{a^{2s + 2 - 2j}}}}} } \left( {\begin{array}{*{20}{c}}
   {2m - 2k + 2j}  \\
   {2j - 1}  \\
\end{array}} \right)\frac{1}{{{n^{2m + 2j - 2k + 1}}}} - \frac{1}{{2{a^{2s + 1}}}}\left\{ {\frac{1}{{{{\left( {n - a} \right)}^{2m + 2}}}} - \frac{1}{{{{\left( {n + a} \right)}^{2m + 2}}}}} \right\}\\
&+ \sum\limits_{j = 1}^s {\frac{1}{{{a^{2s + 2 - 2j}}}}} \left( {\begin{array}{*{20}{c}}
   {2m + 2j}  \\
   {2j - 1}  \\
\end{array}} \right)\frac{1}{{{n^{2m + 2j + 1}}}}.
\end{align*}
At a negative integer $-n$ the pole is simple and residue is
\[\frac{{{\zeta _n}\left( {2m + 1} \right) - \zeta \left( {2m + 1} \right)}}{{{n^{2s}}\left( {{n^2} - {a^2}} \right)}} - \frac{1}{{{n^{2s + 2m + 1}}\left( {{n^2} - {a^2}} \right)}}.\]
The residue of the pole at $\pm a$ is
\[ - \frac{{\pi \cot \left( {\pi a} \right)}}{{2{a^{2s + 1}}}}\sum\limits_{n = 1}^\infty  {\left\{ {\frac{1}{{{{\left( {n - a} \right)}^{2m + 1}}}} + \frac{1}{{{{\left( {n + a} \right)}^{2m + 1}}}}} \right\}} .\]
Finally the residue of the pole of order $2m+2s+2$ at 0 is found to be
\[\sum\limits_{k = 1}^{s + 1} {\left( {\begin{array}{*{20}{c}}
   {2m + 2k - 2}  \\
   {2k - 2}  \\
\end{array}} \right)} \frac{{\zeta \left( {2m + 2k - 1} \right)}}{{{a^{2s - 2k + 4}}}} - 2\sum\limits_{n = 1}^s {\sum\limits_{k = 1}^n {\left( {\begin{array}{*{20}{c}}
   {2m + 2k - 1}  \\
   {2k - 1}  \\
\end{array}} \right)} } \frac{{\zeta \left( {2m + 2k - 1} \right)\zeta \left( {2n - 2k + 2} \right)}}{{{a^{2s - 2n + 2}}}}.\]
Summing these four contributions yields the statement of the theorem.\hfill$\square$\\
Taking $s=0$ in Theorem 2.4, we have the following Corollary.
\begin{cor}([3]) For integers $m\in \N_0$ with $a$ is a real and $a \ne 0,  - 1, - 2, \cdots$, we have
\begin{align*}
\sum\limits_{n = 1}^\infty  {\frac{{{\zeta _n}\left( {2m + 1} \right)}}{{{n^2} - {a^2}}}}
 &=\frac{1}{2}\sum\limits_{n = 1}^\infty  {\frac{1}{{{n^{ 2m + 1}}\left( {{n^2} - {a^2}} \right)}}} \\
 &\quad -\frac{1}{{2a}}\sum\limits_{k = 0}^m {\zeta \left( {2k} \right)} \sum\limits_{n = 1}^\infty  {\left\{ {\frac{1}{{{{\left( {n - a} \right)}^{2m - 2k + 2}}}} - \frac{1}{{{{\left( {n + a} \right)}^{2m - 2k + 2}}}}} \right\}}
\nonumber \\
           &\quad  + \frac{{\pi \cot \left( {\pi a} \right)}}{{4a}}\sum\limits_{n = 1}^\infty  {\left\{ {\frac{1}{{{{\left( {n - a} \right)}^{2m + 1}}}} + \frac{1}{{{{\left( {n + a} \right)}^{2m + 1}}}} - \frac{2}{{{n^{2m + 1}}}}} \right\}}  .\tag{2.16}
\end{align*}
\end{cor}
Note that formula above was also proved in [3] by another method.\\
\indent Next we consider the following type of parametric Euler Sums by the method of constructing function ${T_{s,t}}\left( {x,y} \right)$
\[\sum\limits_{n = 1}^\infty  {\frac{{{H_n}}}{{{{\left( {n + a} \right)}^s}}}} ,\sum\limits_{n = 1}^\infty  {\frac{{{L_n}\left( 1 \right)}}{{{{\left( {n + a} \right)}^s}}}} ,\;a \ne  - 1, - 2, \cdots .\]
\begin{thm} If $a$ is a real number, $s$ is a positive integer and $a \ne  - 1, - 2, \cdots $, then
\begin{align*}
\sum\limits_{n = 1}^\infty  {\frac{{{y^n}\sum\limits_{j = 1}^{n - 1} {\frac{{{x^{n - j}}}}{j}}  + {x^n}\sum\limits_{j = 1}^{n - 1} {\frac{{{y^{n - j}}}}{j}} }}{{{{\left( {n + a} \right)}^s}}}}  =& s{\rm Li}{_{s + 1}}\left( {a,xy} \right) - \sum\limits_{j = 1}^s {{\rm Li}{_j}\left( {a,x} \right){\rm Li}{_{s + 1 - j}}\left( {a,y} \right)}\\&  + {\rm Li}{_s}\left( {a,xy} \right)\left( {{\rm Li}{_1}\left( x \right) + {\rm Li}{_1}\left( y \right)} \right).\tag{2.17}
\end{align*}
where $x,y \in \left[ { - 1,1} \right)$ and the parametric polylogarithm function ${{\rm Li}{_s}\left( {a,x} \right)}$ is defined by
\[{\rm Li}{_s}\left( {a,x} \right) = \sum\limits_{n = 1}^\infty  {\frac{{{x^n}}}{{{{\left( {n + a} \right)}^s}}}} ,\;\Re (s)\geq1, - 1 \le x < 1.\]
If $a=0$, then the  function ${{\rm Li}{_s}\left( {a,x} \right)}$ reduces to the classical polylogarithm function ${{\rm Li}{_s}\left( {x} \right)}$ which is defined by
\[{\rm Li}{_s}\left( x \right) = \sum\limits_{n = 1}^\infty  {\frac{{{x^n}}}{{{n^s}}}} ,\;\Re (s)\geq1, - 1 \le x < 1,\]
with ${\rm Li}_1(x)=-\ln(1-x).$
\end{thm}
\pf Motivated by [2,7], for real $ - 1 \le x,y < 1$ and integers $s$ and $t$, we consider the function
\begin{align*}
{T_{s,t}}\left( {x,y} \right)
 &=\sum\limits_{\scriptstyle m,n = 1 \hfill \atop
  \scriptstyle m \ne n \hfill}^\infty  {\frac{{{x^n}{y^m}}}{{{{\left( {n + a} \right)}^s}{{\left( {m + a} \right)}^t}\left( {m - n} \right)}}}  = \sum\limits_{\scriptstyle m,n = 1 \hfill \atop
  \scriptstyle m \ne n \hfill}^\infty  {\frac{{{x^n}{y^m}\left( {m - n + n + a} \right)}}{{{{\left( {n + a} \right)}^s}{{\left( {m + a} \right)}^{t + 1}}\left( {m - n} \right)}}}
\nonumber \\
& = \sum\limits_{\scriptstyle m,n = 1 \hfill \atop
  \scriptstyle m \ne n \hfill}^\infty  {\frac{{{x^n}{y^m}}}{{{{\left( {n + a} \right)}^s}{{\left( {m + a} \right)}^{t + 1}}}} + \sum\limits_{\scriptstyle m,n = 1 \hfill \atop
  \scriptstyle m \ne n \hfill}^\infty  {\frac{{{x^n}{y^m}}}{{{{\left( {n + a} \right)}^{s - 1}}{{\left( {m + a} \right)}^{t + 1}}\left( {m - n} \right)}}} }
\nonumber \\
& = \sum\limits_{n = 1}^\infty  {\frac{{{x^n}}}{{{{\left( {n + a} \right)}^s}}}} \left\{ {\sum\limits_{m = 1}^\infty  {\frac{{{y^m}}}{{{{\left( {m + a} \right)}^{t + 1}}}}}  - \frac{{{y^n}}}{{{{\left( {n + a} \right)}^{t + 1}}}}} \right\} + {T_{s - 1,t + 1}}\left( {x,y} \right)
\nonumber \\
&= {\rm Li}{_s}\left( {a,x} \right){\rm Li}{_{t + 1}}\left( {a,y} \right) - {\rm Li}{_{s + t + 1}}\left( {a,xy} \right) + {T_{s - 1,t + 1}}\left( {x,y} \right)
.
\end{align*}
Telescoping this gives
\[{T_{s,t}}\left( {x,y} \right) = {T_{0,s + t}}\left( {x,y} \right) - s{\rm Li}{_{s + t + 1}}\left( {a,xy} \right) + \sum\limits_{j = 1}^s {{\rm Li}{_j}\left( {a,x} \right){\rm Li}{_{s+t+1 - j}}\left( {a,y} \right)} ,\;\; s \in \N_0.\]
With $t=0$, this becomes
\[{T_{s,0}}\left( {x,y} \right) = {T_{0,s}}\left( {x,y} \right) - s{\rm Li}{_{s + 1}}\left( {a,xy} \right) + \sum\limits_{j = 1}^s {{\rm Li}{_j}\left( {a,x} \right){\rm Li}{_{s + 1 - j}}\left( {a,y} \right)} ,\;\; s \in \N.\tag{2.18}\]
But for any integers $s$ and $t$, there holds
\begin{align*}
{T_{t,s}}\left( {x,y} \right)& = \sum\limits_{\scriptstyle m,n = 1 \hfill \atop
  \scriptstyle m \ne n \hfill}^\infty  {\frac{{{x^n}{y^m}}}{{{{\left( {n + a} \right)}^t}{{\left( {m + a} \right)}^s}\left( {m - n} \right)}}} \\& =  - \sum\limits_{\scriptstyle m,n = 1 \hfill \atop
  \scriptstyle m \ne n \hfill}^\infty  {\frac{{{y^m}{x^n}}}{{{{\left( {m + a} \right)}^s}{{\left( {n + a} \right)}^t}\left( {n - m} \right)}}}  \\&=  - {T_{s,t}}\left( {y,x} \right).\tag{2.19}
  \end{align*}
From the definition of ${T_{s,t}}\left( {x,y} \right)$, we can deduce that
\begin{align*}
{T_{s,0}}\left( {x,y} \right)
 &=\sum\limits_{\scriptstyle m,n = 1 \hfill \atop
  \scriptstyle m \ne n \hfill}^\infty  {\frac{{{x^n}{y^m}}}{{{{\left( {n + a} \right)}^s}\left( {m - n} \right)}}}  = \sum\limits_{n = 1}^\infty  {\frac{{{x^n}{y^n}}}{{{{\left( {n + a} \right)}^s}}}} \sum\limits_{m = n + 1}^\infty  {\frac{{{y^{m - n}}}}{{m - n}}}  - \sum\limits_{n = 1}^\infty  {\frac{{{x^n}}}{{{{\left( {n + a} \right)}^s}}}} \sum\limits_{m = 1}^{n - 1} {\frac{{{y^m}}}{{n - m}}}
\nonumber \\
&={\rm Li} {_s}\left( {a,xy} \right){\rm Li}{_1}\left( y \right) - \sum\limits_{n = 1}^\infty  {\frac{{{x^n}}}{{{{\left( {n + a} \right)}^s}}}} \sum\limits_{j = 1}^{n - 1} {\frac{{{y^{n - j}}}}{j}}.\ \ \ \ \ \ \ \ \ \ \ \ \ \ \ \ \ \ \ \ \ \ \ \ \ \ \ \ \ \ \ \ \ \ \ \ \ \ \ \ \ \ \ \ \ (2.20)
\end{align*}
Substituting (2.19), (2.20) into (2.18) yields the desired result.\hfill$\square$\\
Setting $x=y$ in (2.17), by simple calculation, we obtain the result
\begin{align*}
\sum\limits_{n = 1}^\infty  {\frac{{{x^n}}}{{{{\left( {n + a} \right)}^s}}}} \sum\limits_{j = 1}^{n - 1} {\frac{{{x^{n - j}}}}{j}}
 &=\frac{s}{2}{\rm Li}{_{s + 1}}\left( {a,{x^2}} \right) + {\rm Li}{_s}\left( {a,{x^2}} \right){\rm Li}{_1}\left( x \right) - {\rm Li}{_s}\left( {a,x} \right){\rm Li}{_1}\left( {a,x} \right)
\nonumber \\
           &\quad - \frac{1}{2}\sum\limits_{j = 2}^{s - 1} {{\rm Li}{_j}\left( {a,x} \right){\rm Li}{_{s + 1 - j}}\left( {a,x} \right)} ,\tag{2.21}
\end{align*}
where in (2.21), we now require $s\in \mathbb{N} \setminus \{1\}$ because the terms $j=1$ and $j=s$ were separated, and assumed to be distinct. Taking $x=-1$ in (2.21), we get
\begin{align*}
\sum\limits_{n = 1}^\infty  {\frac{{{L_n}\left( 1 \right)}}{{{{\left( {n + a} \right)}^s}}}}  = & \frac{1}{2}\sum\limits_{j = 1}^{s - 2} {\bar \zeta \left( {s - j,a + 1} \right)} \bar \zeta \left( {j + 1,a + 1} \right) - \frac{s}{2}\zeta \left( {s + 1,a + 1} \right)\\
& + \zeta \left( {s,a + 1} \right)\ln 2 + \bar \zeta \left( {s,a + 1} \right)\bar \zeta \left( {1,a + 1} \right) + \sum\limits_{n = 1}^\infty  {\frac{{{{\left( { - 1} \right)}^{n - 1}}}}{{n{{\left( {n + a} \right)}^s}}}} .\tag{2.22}
\end{align*}
Note that when $s>1$, $\mathop {\lim }\limits_{x \to 1}{\rm{L}}{{\rm{i}}_s}\left( {a,1} \right) = \zeta \left( {s,a + 1} \right)$; when $s\geq 1$, ${\rm{L}}{{\rm{i}}_s}\left( {a, - 1} \right) =  - \bar \zeta \left( {s,a + 1} \right)$. When $x$ approach 1, we arrive at the conclusion that
\[\mathop {\lim }\limits_{x \to 1} \left\{ {{\rm{L}}{{\rm{i}}_s}\left( {a,{x^2}} \right){\rm{L}}{{\rm{i}}_1}\left( x \right) - {\rm{L}}{{\rm{i}}_s}\left( {a,{x}} \right){\rm{L}}{{\rm{i}}_1}\left( {a,x} \right)} \right\} = a\zeta \left( {s,a + 1} \right)\sum\limits_{n = 1}^\infty  {\frac{1}{{n\left( {n + a} \right)}}} .\tag{2.23}\]
Hence, letting $x\rightarrow 1$ in (2.21), we obtain
\begin{align*}
\sum\limits_{n = 1}^\infty  {\frac{{{H_n}}}{{{{\left( {n + a} \right)}^s}}}}
 &=\frac{s}{2}\zeta \left( {s + 1,a + 1} \right) - \frac{1}{2}\sum\limits_{j = 1}^{s - 2} {\zeta \left( {s - j,a + 1} \right)} \zeta \left( {j + 1,a + 1} \right)
\nonumber \\
           &\quad + a\zeta \left( {s,a + 1} \right)\sum\limits_{n = 1}^\infty  {\frac{1}{{n\left( {n + a} \right)}}}  + \sum\limits_{n = 1}^\infty  {\frac{1}{{n{{\left( {n + a} \right)}^s}}}}.\tag{2.24}
\end{align*}
Putting $a=0$ in (2.22) and (2.24),  we can deduce the well-known identities
\[\sum\limits_{n = 1}^\infty  {\frac{{{L_n}\left( 1 \right)}}{{{n^s}}}}  = \zeta \left( s \right)\ln 2 - \frac{s}{2}\zeta \left( {s + 1} \right) + \bar \zeta \left( {s + 1} \right) + \frac{1}{2}\sum\limits_{j = 1}^s {\bar \zeta \left( {s - j + 1} \right)} \bar \zeta \left( j \right),\]
\[\sum\limits_{n = 1}^\infty  {\frac{{{H_n}}}{{{n^s}}}}  = \frac{1}{2}\left\{ {\left( {s + 2} \right)\zeta \left( {s+ 1} \right) - \sum\limits_{i = 1}^{s - 2} {\zeta \left( {s - i} \right)\zeta \left( {i + 1} \right)} } \right\}.\]
In fact, from (2.21), we can obtain many other evaluation of parametric linear Euler sums. For example,
\begin{align*}
\sum\limits_{n = 1}^\infty  {\frac{{{x^n}}}{{{n^2} - {a^2}}}} \sum\limits_{j = 1}^{n - 1} {\frac{{{x^{n - j}}}}{j}}  &= \frac{1}{{2a}}\left\{ {\sum\limits_{n = 1}^\infty  {\frac{{{x^n}}}{{n - a}}} \sum\limits_{j = 1}^{n - 1} {\frac{{{x^{n - j}}}}{j}}  - \sum\limits_{n = 1}^\infty  {\frac{{{x^n}}}{{n + a}}} \sum\limits_{j = 1}^{n - 1} {\frac{{{x^{n - j}}}}{j}} } \right\}\\
& = \sum\limits_{n = 1}^\infty  {\frac{{n{x^{2n}}}}{{{{\left( {{n^2} - {a^2}} \right)}^2}}}}  + {\rm{L}}{{\rm{i}}_1}\left( x \right)\sum\limits_{n = 1}^\infty  {\frac{{{x^{2n}}}}{{{n^2} - {a^2}}}}  - \left( {\sum\limits_{n = 1}^\infty  {\frac{{n{x^n}}}{{{n^2} - {a^2}}}} } \right)\left( {\sum\limits_{n = 1}^\infty  {\frac{{{x^n}}}{{{n^2} - {a^2}}}} } \right).\tag{2.25}
\end{align*}
By a direct calculation, we can find that the formula (2.25) can be rewritten as
\begin{align*}
\sum\limits_{n = 1}^\infty  {\frac{{{x^n}}}{{{n^2} - {a^2}}}} \sum\limits_{j = 1}^{n - 1} {\frac{{{x^{n - j}}}}{j}}
 &= \sum\limits_{n = 1}^\infty  {\frac{{n{x^{2n}}}}{{{{\left( {{n^2} - {a^2}} \right)}^2}}}}  - a{\left( {\sum\limits_{n = 1}^\infty  {\frac{{{x^n}}}{{{n^2} - {a^2}}}} } \right)^2} + \sum\limits_{n = 1}^\infty  {\frac{{{x^{2n}}{\rm{Li}}{_1}\left( x \right)}}{{{n^2} - {a^2}}}}
\nonumber \\
           &\quad - \left( {\sum\limits_{n = 1}^\infty  {\frac{{{x^n}}}{{n + a}}} } \right)\left( {\sum\limits_{n = 1}^\infty  {\frac{{{x^n}}}{{{n^2} - {a^2}}}} } \right).\tag{2.26}
\end{align*}
The following identity is easily derived
\[\mathop {\lim }\limits_{x \to 1} \left\{ {{\rm{L}}{{\rm{i}}_1}\left( x \right)\sum\limits_{n = 1}^\infty  {\frac{{{x^{2n}}}}{{{n^2} - {a^2}}}}  - \left( {\sum\limits_{n = 1}^\infty  {\frac{{{x^n}}}{{n + a}}} } \right)\left( {\sum\limits_{n = 1}^\infty  {\frac{{{x^n}}}{{{n^2} - {a^2}}}} } \right)} \right\} = a\left( {\sum\limits_{n = 1}^\infty  {\frac{1}{{{n^2} - {a^2}}}} } \right)\left( {\sum\limits_{n = 1}^\infty  {\frac{1}{{n\left( {n + a} \right)}}} } \right).\tag{2.27}\]
Letting $x$ approach 1 in (2.26), we obtain the parametric linear Euler sums
\[\sum\limits_{n = 1}^\infty  {\frac{{{H_n}}}{{{n^2} - {a^2}}}}  = \sum\limits_{n = 1}^\infty  {\frac{n}{{{{\left( {{n^2} - {a^2}} \right)}^2}}}}  + \left( {1 - {a^2}\left( {\sum\limits_{n = 1}^\infty  {\frac{1}{{{n^2} - {a^2}}}} } \right)} \right)\left( {\sum\limits_{n = 1}^\infty  {\frac{1}{{n\left( {{n^2} - {a^2}} \right)}}} } \right).\tag{2.28}\]
Similarly, considering the following limit
\[\mathop {\lim }\limits_{x \to 1} \sum\limits_{n = 1}^\infty  {\left( {\frac{{{x^n}}}{{n + a}} + \frac{{{x^n}}}{{n - a}} - 2\frac{{{x^n}}}{n}} \right)\left( {\sum\limits_{j = 1}^{n - 1} {\frac{{{x^{n - j}}}}{j}} } \right)}  = 2{a^2}\sum\limits_{n = 1}^\infty  {\frac{{{H_{n - 1}}}}{{n\left( {{n^2} - {a^2}} \right)}}},\tag{2.29}\]
and using (2.21), we can deduce the result
\begin{align*}
\sum\limits_{n = 1}^\infty  {\frac{{{H_n}}}{{n\left( {{n^2} - {a^2}} \right)}}}
 &= \frac{3}{2}\sum\limits_{n = 1}^\infty  {\frac{1}{{{{\left( {{n^2} - {a^2}} \right)}^2}}}}  + \sum\limits_{n = 1}^\infty  {\frac{1}{{{n^2}\left( {{n^2} - {a^2}} \right)}}}  - \frac{1}{2}{\left( {\sum\limits_{n = 1}^\infty  {\frac{1}{{{n^2} - {a^2}}}} } \right)^2}
\nonumber \\
           &\quad - \frac{{{a^2}}}{2}\sum\limits_{n = 1}^\infty  {\frac{1}{{{n^2}\left( {{n^2} - {a^2}} \right)^2}}}  - \frac{{{a^2}}}{2}{\left( {\sum\limits_{n = 1}^\infty  {\frac{1}{{{n^2}\left( {{n^2} - {a^2}} \right)}}} } \right)^2}.\tag{2.30}
\end{align*}
\section{Parametric quadratic and cubic Euler sums}
In this section we consider the following type of parametric nonlinear Euler sums involving harmonic numbers by the method of constructing function ${T_{s,t}}\left( {x,y} \right)$
\[\sum\limits_{n = 1}^\infty  {\frac{{H_n^2 - {\zeta _n}\left( 2 \right)}}{{{{\left( {n + a} \right)}^s}}}} ,\;\sum\limits_{n = 1}^\infty  {\frac{{H_n^3 - 3{H_n}{\zeta _n}\left( 2 \right)}}{{{{\left( {n + a} \right)}^s}}}} ,\;s \in \mathbb{N} \setminus \{1\},\;a \ne  - 1, - 2, \cdots .\]
\begin{thm} If $a$ is a real number with $a \ne  - 1, - 2, \cdots $, $s\in \mathbb{N} \setminus \{1\}$. Then
\begin{align*}
\frac{3}{2}\sum\limits_{n = 1}^\infty  {\frac{{H_n^2 - {\zeta _n}\left( 2 \right)}}{{{{\left( {n + a} \right)}^s}}}}
 &=s\sum\limits_{n = 1}^\infty  {\frac{{{H_n}}}{{{{\left( {n + a} \right)}^{s + 1}}}}} +\sum\limits_{n = 1}^\infty  {\frac{{{H_n}}}{{{{n\left( {n + a} \right)}^{s}}}}} - \sum\limits_{j = 2}^{s - 1} {\left( {\sum\limits_{n = 1}^\infty  {\frac{{{H_n}}}{{{{\left( {n + a} \right)}^j}}}} } \right)\zeta \left( {s + 1 - j,a + 1} \right)}
\nonumber \\
           &\quad + a\left( {\sum\limits_{n = 1}^\infty  {\frac{{{H_n}}}{{{{\left( {n + a} \right)}^s}}}} } \right)\left( {\sum\limits_{n = 1}^\infty  {\frac{1}{{n\left( {n + a} \right)}}} } \right) + a\zeta \left( {s,a + 1} \right)\left( {\sum\limits_{n = 1}^\infty  {\frac{{{H_n}}}{{n\left( {n + a} \right)}}} } \right) .\tag{3.1}
\end{align*}
\end{thm}
\pf For real $ - 1 \le x< 1$ and integers $s\ ( s\geqslant 2)$ and $t\ (t\geqslant0)$, consider
\[{T_{s,t}}\left( x \right) := \sum\limits_{\scriptstyle m,n = 1 \hfill \atop
  \scriptstyle m \ne n \hfill}^\infty  {\frac{{{H_n}{x^{n + m}}}}{{{{\left( {n + a} \right)}^s}{{\left( {m + a} \right)}^t}\left( {m - n} \right)}}} .\tag{3.2}\]
Similarly to the proof of Theorem 2.6, we can deduce that
\[{T_{s,t}}\left( x \right) = {T_{0,s + t}}\left( x \right) - sH_{s + t + 1}^{\left( 1 \right)}\left( {a,{x^2}} \right) + \sum\limits_{j = 1}^s {H_j^{\left( 1 \right)}\left( {a,x} \right)H_{s + t + 1 - j}^{\left( 0 \right)}\left( {a,x} \right)}.\tag{3.3} \]
where the function $H_s^{\left( k \right)}\left( {a,x} \right)$ is defined by
\[H_s^{\left( k \right)}\left( {a,x} \right) = \sum\limits_{n = 1}^\infty  {\frac{{H_n^k{x^n}}}{{{{\left( {n + a} \right)}^s}}}} ,\; s\in \N ,\ k \in \N_0,\; - 1 < x < 1.\tag{3.4}\]
Taking $t=0$ in (3.3) yields
\[{T_{s,0}}\left( x \right) = {T_{0,s}}\left( x \right) - sH_{s + 1}^{\left( 1 \right)}\left( {a,{x^2}} \right) + \sum\limits_{j = 1}^s {H_j^{\left( 1 \right)}\left( {a,x} \right)H_{s + 1 - j}^{\left( 0 \right)}\left( {a,x} \right)} .\tag{3.5}\]
By using the definition of ${T_{s,t}}\left( x \right)$, we obtain
\[{T_{s,0}}\left( x \right) = H_s^{\left( 1 \right)}\left( {a,{x^2}} \right){\rm Li}{_1}\left( x \right) - \sum\limits_{n = 1}^\infty  {\frac{{{H_n}{x^n}}}{{{{\left( {n + a} \right)}^s}}}} \left( {\sum\limits_{j = 1}^{n - 1} {\frac{{{x^{n - j}}}}{j}} } \right).\tag{3.6}\]
and
\begin{align*}
{T_{0,s}}\left( x \right) =& \sum\limits_{n = 1}^\infty  {\frac{{{x^n}}}{{{{\left( {n + a} \right)}^s}}}} \left( {\sum\limits_{j = 1}^{n - 1} {\frac{{{H_j}{x^j}}}{{n - j}}} } \right) - H_s^{\left( 0 \right)}\left( {a,{x^2}} \right)H_1^{\left( 1 \right)}\left( {0,x} \right)\\& - \sum\limits_{n = 1}^\infty  {\frac{{{x^{2n}}}}{{{{\left( {n + a} \right)}^s}}}} \left( {\sum\limits_{k = 1}^n {\frac{1}{k}\sum\limits_{j = 1}^\infty  {{x^j}\left( {\frac{1}{j} - \frac{1}{{k + j}}} \right)} } } \right).\tag{3.7}
\end{align*}
Combining (3.5), (3.6) with (3.7), we have
\begin{align*}
&\sum\limits_{n = 1}^\infty  {\frac{{{H_n}{x^n}}}{{{{\left( {n + a} \right)}^s}}}} \left( {\sum\limits_{j = 1}^{n - 1} {\frac{{{x^{n - j}}}}{j}} } \right) + \sum\limits_{n = 1}^\infty  {\frac{{{x^n}}}{{{{\left( {n + a} \right)}^s}}}} \left( {\sum\limits_{j = 1}^{n - 1} {\frac{{{H_j}{x^j}}}{{n - j}}} } \right)\\
& = sH_{s + 1}^{\left( 1 \right)}\left( {a,{x^2}} \right) + \sum\limits_{n = 1}^\infty  {\frac{{{x^{2n}}}}{{{{\left( {n + a} \right)}^s}}}} \left( {\sum\limits_{k = 1}^n {\frac{1}{k}\sum\limits_{j = 1}^\infty  {{x^j}\left( {\frac{1}{j} - \frac{1}{{k + j}}} \right)} } } \right) \\
&\quad- \sum\limits_{j = 2}^{s - 1} {H_j^{\left( 1 \right)}\left( {a,x} \right)H_{s + 1 - j}^{\left( 0 \right)}\left( {a,x} \right)} \\
&\ \ +\left( {H_s^{\left( 1 \right)}\left( {a,{x^2}} \right)L{i_1}\left( x \right) - H_s^{\left( 1 \right)}\left( {a,x} \right)H_1^{\left( 0 \right)}\left( {a,x} \right)} \right)\\
 &\quad+ \left( {H_s^{\left( 0 \right)}\left( {a,{x^2}} \right)H_1^{\left( 1 \right)}\left( {0,x} \right) - H_s^{\left( 0 \right)}\left( {a,x} \right)H_1^{\left( 1 \right)}\left( {a,x} \right)} \right).
\tag{3.8}
\end{align*}
By a direct calculation, we arrive at the conclusion that
\[\mathop {\lim }\limits_{x \to 1} \left( {H_s^{\left( 1 \right)}\left( {a,{x^2}} \right)L{i_1}\left( x \right) - H_s^{\left( 1 \right)}\left( {a,x} \right)H_1^{\left( 0 \right)}\left( {a,x} \right)} \right) = a\left( {\sum\limits_{n = 1}^\infty  {\frac{{{H_n}}}{{{{\left( {n + a} \right)}^s}}}} } \right)\left( {\sum\limits_{n = 1}^\infty  {\frac{1}{{n\left( {n + a} \right)}}} } \right),\tag{3.9}\]
\[\mathop {\lim }\limits_{x \to 1} \left( {H_s^{\left( 0 \right)}\left( {a,{x^2}} \right)H_1^{\left( 1 \right)}\left( {0,x} \right) - H_s^{\left( 0 \right)}\left( {a,x} \right)H_1^{\left( 1 \right)}\left( {a,x} \right)} \right) = a\zeta \left( {s,a + 1} \right)\left( {\sum\limits_{n = 1}^\infty  {\frac{{{H_n}}}{{n\left( {n + a} \right)}}} } \right),\tag{3.10}\]
\[\mathop {\lim }\limits_{x \to 1} \sum\limits_{k = 1}^n {\frac{1}{k}\sum\limits_{j = 1}^\infty  {{x^j}\left( {\frac{1}{j} - \frac{1}{{k + j}}} \right)} }  = \sum\limits_{k = 1}^n {\frac{{{H_k}}}{k}} .\tag{3.11}\]
Taking the limit in (3.8) and combining (2.8), (3.9), (3.10) with (3.11)  yields the desired result. This completes the proof of Theorem 3.1. \hfill$\square$\\
Setting $a=0$ in (3.1), we deduce that well-known identity [1,4,9]
\[\frac{3}{2}\sum\limits_{n = 1}^\infty  {\frac{{H_n^2 - {\zeta _n}\left( 2 \right)}}{{{n^s}}}}  = \left( {s + 1} \right)\sum\limits_{n = 1}^\infty  {\frac{{{H_n}}}{{{n^{s + 1}}}}}  - \sum\limits_{j = 2}^{s - 1} {\left( {\sum\limits_{n = 1}^\infty  {\frac{{{H_n}}}{{{n^j}}}} } \right)\zeta \left( {s + 1 - j} \right)},\;s \in \mathbb{N} \setminus \{1\}.
\tag{3.12}\]
\indent We now evaluate the parametric cubic Euler sums by means of constructing function ${T_{s,t}}\left( {x}\right)$.
\begin{thm} If $a$ is a real number with $a \ne  - 1, - 2, \cdots $, $s \in \mathbb{N} \setminus \{1\} $. Then
\begin{align*}
\sum\limits_{n = 1}^\infty  {\frac{{H_n^3 - 3{H_n}{\zeta _n}\left( 2 \right)}}{{{{\left( {n + a} \right)}^s}}}}
 &=s\sum\limits_{n = 1}^\infty  {\frac{{H_n^2}}{{{{\left( {n + a} \right)}^{s + 1}}}}}  - \sum\limits_{j = 2}^{s - 1} {\left( {\sum\limits_{n = 1}^\infty  {\frac{{{H_n}}}{{{{\left( {n + a} \right)}^j}}}} } \right)\left( {\sum\limits_{n = 1}^\infty  {\frac{{{H_n}}}{{{{\left( {n + a} \right)}^{s + 1 - j}}}}} } \right)}
\nonumber \\
           &\quad  + 2a\left( {\sum\limits_{n = 1}^\infty  {\frac{{{H_n}}}{{{{\left( {n + a} \right)}^s}}}} } \right)\left( {\sum\limits_{n = 1}^\infty  {\frac{{{H_n}}}{{n\left( {n + a} \right)}}} } \right) .\tag{3.13}
\end{align*}
\end{thm}
\pf For real $ - 1 \le x< 1$ and integers $s\ ( s\geqslant 2)$ and $t\ (t\geqslant0)$, we consider
\[{T_{s,t}}\left( x \right) := \sum\limits_{\scriptstyle m,n = 1 \hfill \atop
  \scriptstyle m \ne n \hfill}^\infty  {\frac{{{{H_n}{H_m}}{x^{n + m}}}}{{{{\left( {n + a} \right)}^s}{{\left( {m + a} \right)}^t}\left( {m - n} \right)}}} .\tag{3.14}\]
By a similar argument as inthe proof of Theorem 2.6 and 3.1, we have
\[{T_{s,t}}\left( x \right) = {T_{0,s + t}}\left( x \right) - sH_{s + t + 1}^{\left( 2 \right)}\left( {a,{x^2}} \right) + \sum\limits_{j = 1}^s {H_j^{\left( 1 \right)}\left( {a,x} \right)H_{s + t + 1 - j}^{\left( 1 \right)}\left( {a,x} \right)} .\tag{3.15}\]
Let $t=0$ in (3.15), we obtain
\[{T_{s,0}}\left( x \right) = {T_{0,s}}\left( x \right) - sH_{s + 1}^{\left( 2 \right)}\left( {a,{x^2}} \right) + \sum\limits_{j = 1}^s {H_j^{\left( 1 \right)}\left( {a,x} \right)H_{s + 1 - j}^{\left( 1 \right)}\left( {a,x} \right)} .\tag{3.16}\]
For any integers $s$ and $t$, there holds
\[{T_{t,s}}\left( x \right) = \sum\limits_{\scriptstyle m,n = 1 \hfill \atop
  \scriptstyle m \ne n \hfill}^\infty  {\frac{{{H_n}{H_m}{x^{n + m}}}}{{{{\left( {n + a} \right)}^t}{{\left( {m + a} \right)}^s}\left( {m - n} \right)}}}  =  - \sum\limits_{\scriptstyle m,n = 1 \hfill \atop
  \scriptstyle m \ne n \hfill}^\infty  {\frac{{{H_n}{H_m}{x^{n + m}}}}{{{{\left( {n + a} \right)}^s}{{\left( {m + a} \right)}^t}\left( {m - n} \right)}}}  =  - {T_{s,t}}\left( x \right).\tag{3.17}\]
Therefore, we obtain
\[{T_{s,0}}\left( x \right) = \frac{1}{2}\sum\limits_{j = 1}^s {H_j^{\left( 1 \right)}\left( {a,x} \right)H_{s + 1 - j}^{\left( 1 \right)}\left( {a,x} \right)}  - \frac{s}{2}H_{s + 1}^{\left( 2 \right)}\left( {a,{x^2}} \right).\tag{3.18}\]
On the other hand, from (3.14), we get
\[{T_{s,0}}\left( x \right) = \sum\limits_{n = 1}^\infty  {\frac{{{H_n}{x^{2n}}}}{{{{\left( {n + a} \right)}^s}}}} \left\{ {\sum\limits_{m = n + 1}^\infty  {\frac{{{H_m}{x^{m - n}}}}{{m - n}}} } \right\} + \sum\limits_{n = 1}^\infty  {\frac{{{H_n}{x^{2n}}}}{{{{\left( {n + a} \right)}^s}}}} \left\{ {\sum\limits_{m = 1}^{n - 1} {\frac{{{H_m}{x^m}}}{{m - n}}} } \right\}.\tag{3.19}\]
Noting that
\[\sum\limits_{m = n + 1}^\infty  {\frac{{{H_m}{x^{m - n}}}}{{m - n}}}  = \sum\limits_{j = 1}^\infty  {\frac{{{H_j}}}{j}{x^j} + } \sum\limits_{k = 1}^n {\frac{1}{k}\sum\limits_{j = 1}^\infty  {{x^j}\left( {\frac{1}{j} - \frac{1}{{k + j}}} \right)} }.\tag{3.20} \]
Substituting (3.19) and (3.20) into (3.18) respectively, we get
\begin{align*}
\sum\limits_{n = 1}^\infty  {\frac{{{H_n}{x^{2n}}}}{{{{\left( {n + a} \right)}^s}}}} \left\{ {\sum\limits_{m = 1}^{n - 1} {\frac{{{H_m}{x^m}}}{{n - m}}} } \right\}
 &=\frac{s}{2}H_{s + 1}^{\left( 2 \right)}\left( {a,{x^2}} \right) + \sum\limits_{n = 1}^\infty  {\frac{{{H_n}{x^{2n}}}}{{{{\left( {n + a} \right)}^s}}}} \left\{ {\sum\limits_{k = 1}^n {\frac{1}{k}\sum\limits_{j = 1}^\infty  {{x^j}\left( {\frac{1}{j} - \frac{1}{{k + j}}} \right)} } } \right\}
\nonumber \\
           &\quad \ \  + \left( {H_s^{\left( 1 \right)}\left( {a,{x^2}} \right)H_1^{\left( 1 \right)}\left( {0,x} \right) - H_s^{\left( 1 \right)}\left( {a,x} \right)H_1^{\left( 1 \right)}\left( {a,x} \right)} \right)
\nonumber \\
           &\quad \ \
           - \frac{1}{2}\sum\limits_{j = 2}^{s - 1} {H_j^{\left( 1 \right)}\left( {a,x} \right)H_{s + 1 - j}^{\left( 1 \right)}\left( {a,x} \right)}.\tag{3.21}
\end{align*}
A direct calculation yields
\[\mathop {\lim }\limits_{x \to 1} \left( {H_s^{\left( 1 \right)}\left( {a,{x^2}} \right)H_1^{\left( 1 \right)}\left( {0,x} \right) - H_s^{\left( 1 \right)}\left( {a,x} \right)H_1^{\left( 1 \right)}\left( {a,x} \right)} \right) = a\left( {\sum\limits_{n = 1}^\infty  {\frac{{{H_n}}}{{{{\left( {n + a} \right)}^s}}}} } \right)\left( {\sum\limits_{n = 1}^\infty  {\frac{{{H_n}}}{{n\left( {n + a} \right)}}} } \right).\tag{3.22}\]
Letting $x\rightarrow 1$ in (3.21) and using (2.8), (3.11), (3.22) yields the desired result. \hfill$\square$\\
Putting $a=0$ in (3.13), we obtain the result
\[\sum\limits_{n = 1}^\infty  {\frac{{H_n^3 - 3{H_n}{\zeta _n}\left( 2 \right)}}{{{n^s}}}}  = s\sum\limits_{n = 1}^\infty  {\frac{{H_n^2}}{{{n^{s + 1}}}}}  - \sum\limits_{j = 2}^{s - 1} {\left( {\sum\limits_{n = 1}^\infty  {\frac{{{H_n}}}{{{n^j}}}} } \right)\left( {\sum\limits_{n = 1}^\infty  {\frac{{{H_n}}}{{{n^{s + 1 - j}}}}} } \right)} ,\;s \in \mathbb{N} \setminus \{1\}.
\tag{3.23}\]

\section{Other parametric Euler sums}
In this section we consider the following type of linear parametric Euler sums
\[\sum\limits_{n = 1}^\infty  {\frac{1}{{{{\left( {n + b} \right)}^s}}}} \left( {\sum\limits_{k = 1}^n {\frac{{{x^k}}}{{k + a}}} } \right),\; s \in \mathbb{N} \setminus \{1\},\;a,b \ne  - 1, - 2, \cdots .\tag{4.1}\]
with $x \in \left[ { - 1,1} \right]$. Moreover, we find interesting representations for linear, quadratic parametric Euler sums. First, we define the function ${H_m}\left( x,a \right)$ by
\[{H_m}\left( x,a \right) = \sum\limits_{n = 1}^\infty  {\frac{{{x^{n + a}}}}{{{{\left( {n + a} \right)}^m}}}} ,\;m > 1,\;a\neq -1,-2,\cdots,\;x \in \left[ { - 1,1} \right].\tag{4.2}\]
It is obvious that the function ${H_m}\left( x,a \right)$ have the following propositions,
\[\frac{d}{{dx}}\left( {{H_m}\left( x,a \right)} \right) = \frac{1}{x}{H_{m - 1}}\left( x ,a\right),\;{H_m}\left( x,a \right) = {x^a}{\rm Li}{_m}\left( {a,x} \right).\tag{4.3}\]
\[\mathop {\lim }\limits_{x \to 1} {H_m}\left( x,a \right) = \mathop {\lim }\limits_{x \to 1} \sum\limits_{n = 1}^\infty  {\frac{{{x^{n + a}}}}{{{{\left( {n + a} \right)}^m}}}}  = \sum\limits_{n = 1}^\infty  {\frac{1}{{{{\left( {n + a} \right)}^m}}}}  = \zeta \left( {m,a + 1} \right),\;{\mathop{\Re}\nolimits} \left( m \right) > 1.\]
\begin{thm} For integer $m,n\in \N$ and $\left| x \right| < 1$ with $a,b\neq -1,-2,\ldots$,. Then
\begin{align*}
\int\limits_0^x {{H_m}\left( {t,a} \right)} {t^{n + b - 1}}dt =& \sum\limits_{k = 1}^{m - 1} {\frac{{{{\left( { - 1} \right)}^{k - 1}}{x^{n + b}}}}{{{{\left( {n + b} \right)}^k}}}} {H_{m + 1 - k}}\left( {x,a} \right)\\
& + \frac{{{{\left( { - 1} \right)}^{m - 1}}}}{{{{\left( {n + b} \right)}^m}}}\left\{ {{x^{n + b}}{H_1}\left( {x,a} \right) + \sum\limits_{k = 1}^n {\frac{{{x^{k + a + b}}}}{{k + a + b}}}  - {H_1}\left( {x,a + b} \right)} \right\}.\tag{4.4}
\end{align*}
\end{thm}
\pf Using integration by parts, we have the following recurrence relation
\begin{align*}
\int\limits_0^x {{H_m}\left( {t,a} \right)} {t^{n + b - 1}}dt =& \frac{{{x^{n + b}}}}{{n + b}}{H_m}\left( {x,a} \right) - \frac{1}{{n + b}}\int\limits_0^x {{H_{m - 1}}\left( {t,a} \right)} {t^{n + b - 1}}dt=\cdots\\
 =& \sum\limits_{k = 1}^{m - 1} {\frac{{{{\left( { - 1} \right)}^{k - 1}}{x^{n + b}}}}{{{{\left( {n + b} \right)}^k}}}} {H_{m + 1 - k}}\left( {x,a} \right) + \frac{{{{\left( { - 1} \right)}^{m - 1}}}}{{{{\left( {n + b} \right)}^{m - 1}}}}\int\limits_0^x {{H_1}\left( {t,a} \right)} {t^{n + b - 1}}dt.\tag{4.5}
 \end{align*}
By simple calculation, the second term integral on the right hand side is equal to
\[\int\limits_0^x {{H_1}\left( {t,a} \right)} {t^{n + b - 1}}dt = \frac{1}{{n + b}}\left\{ {{x^{n + b}}{H_1}\left( {x,a} \right) + \sum\limits_{k = 1}^n {\frac{{{x^{k + a + b}}}}{{k + a + b}}}  - {H_1}\left( {x,a + b} \right)} \right\}.\tag{4.6}\]
Substituting (4.6) into (4.5) respectively, we may easily deduce the desired result.  \hfill$\square$
\begin{thm} Let $m, p$ be positive integers with $a,b \ne  - 1, - 2, \cdots $ and $x \in \left[ { - 1,1} \right]$, then
\begin{align*}
&\sum\limits_{n = 1}^\infty  {\left\{ {\frac{{{{\left( { - 1} \right)}^{p - 1}}}}{{{{\left( {n + a} \right)}^{m + p}}}} - \frac{{{{\left( { - 1} \right)}^{m - 1}}}}{{{{\left( {n + b} \right)}^{m + p}}}}} \right\}} \left( {\sum\limits_{k = 1}^n {\frac{{{x^{k + a + b}}}}{{k + a + b}}} } \right)\\
& = \sum\limits_{k = 1}^{m - 1} {{{\left( { - 1} \right)}^{k - 1}}} {H_{m + 1 - k}}\left( x ,a\right){H_{p + k}}\left( x ,b\right) - \sum\limits_{k = 1}^{p - 1} {{{\left( { - 1} \right)}^{k - 1}}} {H_{p + 1 - k}}\left( x,b \right){H_{m + k}}\left( x,a \right) \\
&\ \  + {\left( { - 1} \right)^{m - 1}}\left( {{H_{p + m}}\left( x,b \right){H_1}\left( x,a \right) - {H_{p + m}}\left( 1,b \right){H_1}\left( x,a+b \right)} \right)\\
&\ \  - {\left( { - 1} \right)^{p - 1}}\left( {{H_{p + m}}\left( x,a \right){H_1}\left( x,b \right) - {H_{p + m}}\left( 1,a \right){H_1}\left( x,a+b \right)} \right).
\tag{4.7}
\end{align*}
\end{thm}
\pf We consider the following integral
\[\int\limits_0^x {\frac{{{H_m}\left( t,a \right){H_p}\left( t ,b\right)}}{t}} dt = \sum\limits_{n = 1}^\infty  {\frac{1}{{{{\left( {n + a} \right)}^m}}}} \int\limits_0^x {{H_p}\left( t,b \right)} {t^{n + a - 1}}dt = \sum\limits_{n = 1}^\infty  {\frac{1}{{{{\left( {n + b} \right)}^p}}}} \int\limits_0^x {{H_m}\left( t,a \right)} {t^{n + b - 1}}dt.\tag{4.8}\]
Substituting (4.4) into (4.8), we can obtain (4.7).\hfill$\square$\\
Letting $x \to 1$ in (4.7), it is easily show that
\begin{align*}
&\sum\limits_{n = 1}^\infty  {\left\{ {\frac{{{{\left( { - 1} \right)}^{p - 1}}}}{{{{\left( {n + a} \right)}^{m + p}}}} - \frac{{{{\left( { - 1} \right)}^{m - 1}}}}{{{{\left( {n + b} \right)}^{m + p}}}}} \right\}} \left( {\sum\limits_{k = 1}^n {\frac{1}{{k + a + b}}} } \right)\\
& = \sum\limits_{k = 1}^{m - 1} {{{\left( { - 1} \right)}^{k - 1}}} \zeta \left( {m + 1 - k,a + 1} \right)\zeta \left( {p + k,b + 1} \right) \\
&\quad- \sum\limits_{k = 1}^{p - 1} {{{\left( { - 1} \right)}^{k - 1}}} \zeta \left( {p + 1 - k,b + 1} \right)\zeta \left( {m + k,a + 1} \right) \\
&\quad+ b{\left( { - 1} \right)^{m - 1}}\zeta \left( {m + p,b + 1} \right)\left( {\sum\limits_{n = 1}^\infty  {\frac{1}{{\left( {n + a} \right)\left( {n + a + b} \right)}}} } \right)
 \\&\quad- a{\left( { - 1} \right)^{p - 1}}\zeta \left( {m + p,a + 1} \right)\left( {\sum\limits_{n = 1}^\infty  {\frac{1}{{\left( {n + b} \right)\left( {n + a + b} \right)}}} } \right).
\tag{4.9}
\end{align*}
Taking $m=p=1, b=-a$ in (4.7), we obtain the parametric linear Euler sums
\[\sum\limits_{n = 1}^\infty  {\frac{{n{H_n}}}{{{{\left( {{n^2} - {a^2}} \right)}^2}}}}  = \frac{1}{4}\zeta \left( {2,1 + a} \right)\sum\limits_{n = 1}^\infty  {\frac{1}{{n\left( {n - a} \right)}}}  + \frac{1}{4}\zeta \left( {2,1 - a} \right)\sum\limits_{n = 1}^\infty  {\frac{1}{{n\left( {n + a} \right)}}} , a\neq \pm 1,\pm2,\ldots.\tag{4.10}\]
On the other hand, from (2.21), we can deduce that
\begin{align*}
\sum\limits_{n = 1}^\infty  {\frac{{n{x^n}}}{{{{\left( {{n^2} - {a^2}} \right)}^2}}}} \sum\limits_{j = 1}^{n - 1} {\frac{{{x^{n - j}}}}{j}}
 &=\frac{1}{{4a}}\left\{ {{\rm Li}{_2}\left( { - a,{x^2}} \right){\rm Li}{_1}\left( x \right) - {\rm Li}{_2}\left( { - a,x} \right){\rm Li}{_1}\left( { - a,x} \right)} \right\}
\nonumber \\
           &\quad - \frac{1}{{4a}}\left\{ {{\rm Li}{_2}\left( {a,{x^2}} \right){\rm Li}{_1}\left( x \right) - {\rm Li}{_2}\left( {a,x} \right){\rm Li}{_1}\left( {a,x} \right)} \right\}
\nonumber \\
           &\quad +\frac{1}{2}\sum\limits_{n = 1}^\infty  {\frac{{3{n^2} + {a^2}}}{{{{\left( {{n^2} - {a^2}} \right)}^3}}}} {x^{2n}} .\tag{4.11}
\end{align*}
Furthermore, letting $x$ approach 1 in (4.11) and combining (2.23), we obtain
\begin{align*}
\sum\limits_{n = 1}^\infty  {\frac{{n{H_n}}}{{{{\left( {{n^2} - {a^2}} \right)}^2}}}}
 &= \frac{1}{2}\sum\limits_{n = 1}^\infty  {\frac{{3{n^2} + {a^2}}}{{{{\left( {{n^2} - {a^2}} \right)}^3}}}}  - \frac{1}{4}\zeta \left( {2,1 - a} \right)\sum\limits_{n = 1}^\infty  {\frac{1}{{n\left( {n - a} \right)}}}
\nonumber \\
           &\quad + \sum\limits_{n = 1}^\infty  {\frac{1}{{{{\left( {{n^2} - {a^2}} \right)}^2}}}}  - \frac{1}{4}\zeta \left( {2,1 + a} \right)\sum\limits_{n = 1}^\infty  {\frac{1}{{n\left( {n + a} \right)}}} .\tag{4.12}
\end{align*}
Comparing (4.8) and (4.11), we derive the following beautiful result
\[\frac{5}{2}\sum\limits_{n = 1}^\infty  {\frac{1}{{{{\left( {{n^2} - {a^2}} \right)}^2}}}}  - {\left( {\sum\limits_{n = 1}^\infty  {\frac{1}{{{n^2} - {a^2}}}} } \right)^2} = 2{a^2}\left\{ {\sum\limits_{n = 1}^\infty  {\frac{1}{{{n^2} - {a^2}}}} \sum\limits_{n = 1}^\infty  {\frac{1}{{{{\left( {{n^2} - {a^2}} \right)}^2}}}}  - \sum\limits_{n = 1}^\infty  {\frac{1}{{{{\left( {{n^2} - {a^2}} \right)}^3}}}} } \right\}.\tag{4.13}\]
Setting $a=0$ in (4.13), we can obtain the well-known result
$$\zeta^2(2)=\frac 5{2}\zeta(4).$$
From Theorem 4.2, taking $b=a, p=m+1$, we can give the following Corollary.
\begin{cor} For integer $m\in \N$ and $x\in [-1,1]$ with $a\neq -1,-2,\ldots$ . Then
\begin{align*}
\sum\limits_{n = 1}^\infty  {\frac{1}{{{{\left( {n + a} \right)}^{2m + 1}}}}} \left( {\sum\limits_{k = 1}^n {\frac{{{x^k}}}{{k + 2a}}} } \right)
 &=\sum\limits_{k = 1}^{m - 1} {{{\left( { - 1} \right)}^{m + k - 1}}} {\rm Li}{_{m + 1 - k}}\left( {a,x} \right){\rm Li}{_{m + 1 + k}}\left( {a,x} \right)
\nonumber \\
           &\quad - \left( {{\rm Li}{_{2m + 1}}\left( {a,x} \right){\rm Li}{_1}\left( {a,x} \right) - \zeta \left( {2m + 1,a + 1} \right){\rm Li}{_1}\left( {2a,x} \right)} \right)
\nonumber \\
           &\quad+ \frac{{{{\left( { - 1} \right)}^{m - 1}}}}{2}{\rm Li}_{m + 1}^2\left( {a,x} \right) .\tag{4.14}
\end{align*}
\end{cor}
Putting $x=1$ in (4.14), we have
\begin{align*}
\sum\limits_{n = 1}^\infty  {\frac{1}{{{{\left( {n + a} \right)}^{2m + 1}}}}} \left( {\sum\limits_{k = 1}^n {\frac{1}{{k + 2a}}} } \right)
 &=\sum\limits_{k = 1}^{m - 1} {{{\left( { - 1} \right)}^{m + k - 1}}} \zeta \left( {m + 1 - k,a + 1} \right)\zeta \left( {m + 1 + k,a + 1} \right)
\nonumber \\
           &\quad - a\zeta \left( {2m + 1,a + 1} \right)\sum\limits_{n = 1}^\infty  {\frac{1}{{\left( {n + a} \right)\left( {n + 2a} \right)}}}
\nonumber \\
           &\quad+ \frac{{{{\left( { - 1} \right)}^{m - 1}}}}{2}{\zeta ^2}\left( {m + 1,a + 1} \right) .\tag{4.15}
\end{align*}
Next, we establish some connection between linear and quadratic parametric Euler sums.
\begin{thm} Let $m\geq 0, p>1$ be integers with $a,b \ne  - 1, - 2, \cdots $ and $x \in \left[ { - 1,1} \right)$, we have
\begin{align*}
& {\left( { - 1} \right)^{p - 1}}\sum\limits_{n = 1}^\infty  {\left\{ {\frac{{{\zeta _n}\left( {p + 2m,b + 1} \right)}}{{{{\left( {n + b} \right)}^p}}} - \frac{{{\zeta _n}\left( {p,a + 1} \right)}}{{{{\left( {n + a} \right)}^{p + 2m}}}}} \right\}} \left( {\sum\limits_{k = 1}^n {\frac{{{x^{k + a + b}}}}{{k + a + b}}} } \right) \\
  &= \sum\limits_{i = 1}^{p + 2m - 1} {{{\left( { - 1} \right)}^{i - 1}}{H_{p + 2m + 1 - i}}\left( {x,b} \right)\left( {\sum\limits_{n = 1}^\infty  {\frac{{{\zeta _n}\left( {p,a + 1} \right)}}{{{{\left( {n + a} \right)}^i}}}{x^{n + a}}} } \right)}  \\
 & \quad- \sum\limits_{i = 1}^{p - 1} {{{\left( { - 1} \right)}^{i - 1}}{H_{p + 1 - i}}\left( {x,a} \right)\left( {\sum\limits_{n = 1}^\infty  {\frac{{{\zeta _n}\left( {p + 2m,b + 1} \right)}}{{{{\left( {n + b} \right)}^i}}}{x^{n + b}}} } \right)}  \\
  &\quad+ {\left( { - 1} \right)^{p - 1}}\sum\limits_{n = 1}^\infty  {\frac{{{\zeta _n}\left( {p,a + 1} \right)}}{{{{\left( {n + a} \right)}^{p + 2m}}}}\left\{ {{x^{n + a}}{H_1}\left( {x,b} \right) - {H_1}\left( {x,a + b} \right)} \right\}}  \\
  &\quad- {\left( { - 1} \right)^{p - 1}}\sum\limits_{n = 1}^\infty  {\frac{{{\zeta _n}\left( {p + 2m,b + 1} \right)}}{{{{\left( {n + b} \right)}^p}}}\left\{ {{x^{n + b}}{H_1}\left( {x,a} \right) - {H_1}\left( {x,a + b} \right)} \right\}}
,\tag{4.16}
\end{align*}
and
\begin{align*}
 &{\left( { - 1} \right)^{p - 1}}\sum\limits_{n = 1}^\infty  {\left\{ {\frac{{{\zeta _n}\left( {p + 2m + 1,b + 1} \right)}}{{{{\left( {n + b} \right)}^p}}} + \frac{{{\zeta _n}\left( {p,a + 1} \right)}}{{{{\left( {n + a} \right)}^{p + 2m + 1}}}}} \right\}} \left( {\sum\limits_{k = 1}^n {\frac{{{x^{k + a + b}}}}{{k + a + b}}} } \right) \\
  &= \sum\limits_{i = 1}^{p + 2m} {{{\left( { - 1} \right)}^{i - 1}}{H_{p + 2m + 2 - i}}\left( {x,b} \right)\left( {\sum\limits_{n = 1}^\infty  {\frac{{{\zeta _n}\left( {p,a + 1} \right)}}{{{{\left( {n + a} \right)}^i}}}{x^{n + a}}} } \right)}  \\
  &\quad- \sum\limits_{i = 1}^{p - 1} {{{\left( { - 1} \right)}^{i - 1}}{H_{p + 1 - i}}\left( {x,a} \right)\left( {\sum\limits_{n = 1}^\infty  {\frac{{{\zeta _n}\left( {p + 2m + 1,b + 1} \right)}}{{{{\left( {n + b} \right)}^i}}}{x^{n + b}}} } \right)}  \\
 &\quad + {\left( { - 1} \right)^p}\sum\limits_{n = 1}^\infty  {\frac{{{\zeta _n}\left( {p,a+1} \right)}}{{{{\left( {n + a} \right)}^{p + 2m + 1}}}}\left\{ {{x^{n + a}}{H_1}\left( {x,b} \right) - {H_1}\left( {x,a + b} \right)} \right\}}  \\
 &\quad + {\left( { - 1} \right)^p}\sum\limits_{n = 1}^\infty  {\frac{{{\zeta _n}\left( {p + 2m + 1,b+1} \right)}}{{{{\left( {n + b} \right)}^p}}}\left\{ {{x^{n + b}}{H_1}\left( {x,a} \right) - {H_1}\left( {x,a + b} \right)} \right\}}
.\tag{4.17}
\end{align*}
\end{thm}
\pf To prove the first identity, we consider the following integral
\[\int\limits_0^x {\frac{{{H_p}\left( {t,a} \right){H_{p + 2m}}\left( {t,b} \right)}}{{t\left( {1 - t} \right)}}dt},x\in (-1,1). \]
By using the definition of ${H_m}\left( x,a \right)$ and the cauchy product of power series, we have
\[\frac{{{H_p}\left( {t,a} \right)}}{{1 - t}} = \sum\limits_{n = 1}^\infty  {{\zeta _n}\left( {p,a+1} \right){t^{n + a}}} ,\frac{{{H_{p + 2m}}\left( {t,b} \right)}}{{1 - t}} = \sum\limits_{n = 1}^\infty  {{\zeta _n}\left( {p + 2m,b+1} \right){t^{n + b}}} ,\;t \in \left( { - 1,1} \right).\]
Furthermore, by (4.4), we conclude that
\begin{align*}
\int\limits_0^x {\frac{{{H_p}\left( {t,a} \right){H_{p + 2m}}\left( {t,b} \right)}}{{t\left( {1 - t} \right)}}dt}  =& \sum\limits_{n = 1}^\infty  {{\zeta _n}\left( {p,a+1} \right)\int\limits_0^x {{t^{n + a - 1}}{H_{p + 2m}}\left( {t,b} \right)dt} } \\
= &\sum\limits_{i = 1}^{p + 2m - 1} {{{\left( { - 1} \right)}^{i - 1}}{H_{p + 2m + 1 - i}}\left( {x,b} \right)\left( {\sum\limits_{n = 1}^\infty  {\frac{{{\zeta _n}\left( {p,a + 1} \right)}}{{{{\left( {n + a} \right)}^i}}}{x^{n + a}}} } \right)}\\
& + {\left( { - 1} \right)^{p - 1}}\sum\limits_{n = 1}^\infty  {\frac{{{\zeta _n}\left( {p,a + 1} \right)}}{{{{\left( {n + a} \right)}^{p + 2m}}}}\left\{ {{x^{n + a}}{H_1}\left( {x,b} \right) - {H_1}\left( {x,a + b} \right)} \right\}} \\
& + {\left( { - 1} \right)^{p - 1}}\sum\limits_{n = 1}^\infty  {\frac{{{\zeta _n}\left( {p,a + 1} \right)}}{{{{\left( {n + a} \right)}^{p + 2m}}}}} \left( {\sum\limits_{k = 1}^n {\frac{{{x^{k + a + b}}}}{{k + a + b}}} } \right)\\
=&  \sum\limits_{n = 1}^\infty  {{\zeta _n}\left( {p + 2m,b+1} \right)\int\limits_0^x {{t^{n + b - 1}}{H_p}\left( {t,a} \right)dt} } \\
 =& \sum\limits_{i = 1}^{p - 1} {{{\left( { - 1} \right)}^{i - 1}}{H_{p + 1 - i}}\left( {x,a} \right)\left( {\sum\limits_{n = 1}^\infty  {\frac{{{\zeta _n}\left( {p + 2m,b + 1} \right)}}{{{{\left( {n + b} \right)}^i}}}{x^{n + b}}} } \right)} \\
 & + {\left( { - 1} \right)^{p - 1}}\sum\limits_{n = 1}^\infty  {\frac{{{\zeta _n}\left( {p + 2m,b + 1} \right)}}{{{{\left( {n + b} \right)}^p}}}\left\{ {{x^{n + b}}{H_1}\left( {x,a} \right) - {H_1}\left( {x,a + b} \right)} \right\}} \\
 &+ {\left( { - 1} \right)^{p - 1}}\sum\limits_{n = 1}^\infty  {\frac{{{\zeta _n}\left( {p + 2m,b + 1} \right)}}{{{{\left( {n + b} \right)}^p}}}} \left( {\sum\limits_{k = 1}^n {\frac{{{x^{k + a + b}}}}{{k + a + b}}} } \right).
\end{align*}
By a direct calculation, we can deduce the desired result. To prove the second identity of our
theorem, we use the following integral
\[\int\limits_0^x {\frac{{{H_p}\left( {t,a} \right){H_{p + 2m + 1}}\left( {t,b} \right)}}{{t\left( {1 - t} \right)}}dt} ,x\in (-1,1).\]
and apply the same arguments as in the proof of (4.16), we may easily deduce the result.\hfill$\square$\\
We now consider the following function
\[y = \sum\limits_{n = 1}^\infty  {\left\{ {{\zeta _n}\left( {1,a+1} \right){\zeta _n}\left( {p,a+1} \right) - {\zeta _n}\left( {p + 1,a+1} \right)} \right\}{x^{n + a - 1}}} ,\;x \in \left( { - 1,1} \right).\]
From the definition of ${{\zeta _n}\left( {p,a+1} \right)}$, we can find that
\[\left( {1 - x} \right)y = \sum\limits_{n = 1}^\infty  {\left\{ {\frac{{{\zeta _n}\left( {1,a+1} \right)}}{{{{\left( {n + a} \right)}^p}}} + \frac{{{\zeta _n}\left( {p,a+1} \right)}}{{n + a}} - 2\frac{1}{{{{\left( {n + a} \right)}^{p + 1}}}}} \right\}{x^{n + a - 1}}}.\tag{4.18}\]
Using the elementary integral identity
\[\int\limits_0^1 {\frac{{{x^{n + a - 1}}{{\ln }^m}x}}{{1 - x}}} dx = {\left( { - 1} \right)^m}m!\left( {\zeta \left( {m + 1,a + 1} \right) - {\zeta _{n - 1}}\left( {m + 1,a+1} \right)} \right),\: m \in \N,\]
then multiplying (4.18) by ${\frac{{{{\ln }^{m - 1}}x}}{{1 - x}}}$, and integrating over $(0,1)$. The result is
\begin{align*}
 &\sum\limits_{n = 1}^\infty  {\frac{{{\zeta _n}\left( {1,a+1} \right){\zeta _n}\left( {p,a+1} \right) - {\zeta _n}\left( {p + 1,a+1} \right)}}{{{{\left( {n + a} \right)}^m}}}}  \\
  &= \sum\limits_{n = 1}^\infty  {\left\{ {\frac{{{\zeta _n}\left( {1,a+1} \right)}}{{{{\left( {n + a} \right)}^p}}} + \frac{{{\zeta _n}\left( {p,a+1} \right)}}{{n + a}} - 2\frac{1}{{{{\left( {n + a} \right)}^{p + 1}}}}} \right\}\left( {\zeta \left( {m,a + 1} \right) - {\zeta _{n - 1}}\left( {m,a+1} \right)} \right)} .
\end{align*}
By simple calculation, we arrive at the conclusion that
\begin{align*}
&\sum\limits_{n = 1}^\infty  {\left\{ {\frac{{{\zeta _n}\left( {1,a+1} \right){\zeta _n}\left( {p,a+1} \right)}}{{{{\left( {n + a} \right)}^m}}} + \frac{{{\zeta _n}\left( {1,a+1} \right){\zeta _n}\left( {m,a+1} \right)}}{{{{\left( {n + a} \right)}^p}}}} \right\}} \\
& = \zeta \left( {m,a + 1} \right)\sum\limits_{n = 1}^\infty  {\frac{{{\zeta _n}\left( {1,a+1} \right)}}{{{{\left( {n + a} \right)}^p}}}}  + \sum\limits_{n = 1}^\infty  {\frac{{{\zeta _n}\left( {1,a+1} \right)}}{{{{\left( {n + a} \right)}^{p + m}}}}}  + \sum\limits_{n = 1}^\infty  {\frac{{{\zeta _n}\left( {p,a+1} \right)}}{{{{\left( {n + a} \right)}^{m + 1}}}}} \\
&\quad - \sum\limits_{n = 1}^\infty  {\frac{{{\zeta _n}\left( {p + 1,a+1} \right)}}{{{{\left( {n + a} \right)}^m}}}}  + \sum\limits_{n = 1}^\infty  {\frac{{{\zeta _n}\left( {p,a+1} \right)}}{{n + a}}\left( {\zeta \left( {m,a + 1} \right) - {\zeta _n}\left( {m,a+1} \right)} \right)}.\tag{4.19}
\end{align*}
Change $m,p$ to $p,m$, we obtain
\begin{align*}
&\sum\limits_{n = 1}^\infty  {\left\{ {\frac{{{\zeta _n}\left( {1,a+1} \right){\zeta _n}\left( {m,a+1} \right)}}{{{{\left( {n + a} \right)}^p}}} + \frac{{{\zeta _n}\left( {1,a+1} \right){\zeta _n}\left( {p,a+1} \right)}}{{{{\left( {n + a} \right)}^m}}}} \right\}} \\
& = \zeta \left( {p,a + 1} \right)\sum\limits_{n = 1}^\infty  {\frac{{{\zeta _n}\left( {1,a+1} \right)}}{{{{\left( {n + a} \right)}^m}}}}  + \sum\limits_{n = 1}^\infty  {\frac{{{\zeta _n}\left( {1,a+1} \right)}}{{{{\left( {n + a} \right)}^{m + p}}}}}  + \sum\limits_{n = 1}^\infty  {\frac{{{\zeta _n}\left( {m,a+1} \right)}}{{{{\left( {n + a} \right)}^{p + 1}}}}} \\
&\quad - \sum\limits_{n = 1}^\infty  {\frac{{{\zeta _n}\left( {m + 1,a+1} \right)}}{{{{\left( {n + a} \right)}^p}}}}  + \sum\limits_{n = 1}^\infty  {\frac{{{\zeta _n}\left( {m,a+1} \right)}}{{n + a}}\left( {\zeta \left( {p,a + 1} \right) - {\zeta _n}\left( {p,a+1} \right)} \right)} .\tag{4.20}
\end{align*}
Combining (4.19) with (4.20), we have the result
\begin{align*}
&\sum\limits_{n = 1}^\infty  {\frac{{\zeta \left( {m,a + 1} \right){\zeta _n}\left( {p,a+1} \right) - \zeta \left( {p,a + 1} \right){\zeta _n}\left( {m,a+1} \right)}}{{n + a}}} \\
& = \zeta \left( {p,a + 1} \right)\sum\limits_{n = 1}^\infty  {\frac{{{\zeta _n}\left( {1,a+1} \right)}}{{{{\left( {n + a} \right)}^m}}}}  - \zeta \left( {m,a + 1} \right)\sum\limits_{n = 1}^\infty  {\frac{{{\zeta _n}\left( {1,a+1} \right)}}{{{{\left( {n + a} \right)}^p}}}} \\
&\quad + \zeta \left( {m,a + 1} \right)\zeta \left( {p + 1,a + 1} \right) - \zeta \left( {m + 1,a + 1} \right)\zeta \left( {p,a + 1} \right).\tag{4.21}
\end{align*}
Taking $x\rightarrow 1,\ b=a$ in (4.16), (4.17) and using (4.21) with the following identity
\[\mathop {\lim }\limits_{x \to 1} \left\{ {{x^{n + a}}{H_1}\left( x,b \right) - {H_1}\left( x,a+b \right)} \right\} = a\sum\limits_{n = 1}^\infty  {\frac{1}{{\left( {n + b} \right)\left( {n + a + b} \right)}}},\]
we obtain interesting representations for linear, quadratic parametric Euler sums:
\begin{cor}
For integers $p\geq 2,m\geq 0$ and $a\neq -1,-2,\ldots$, we have
\begin{align*}
&{\left( { - 1} \right)^{p - 1}}\sum\limits_{n = 1}^\infty  {\left\{ {\frac{{{\zeta _n}\left( {p + 2m,a + 1} \right)}}{{{{\left( {n + a} \right)}^p}}} - \frac{{{\zeta _n}\left( {p,a + 1} \right)}}{{{{\left( {n + a} \right)}^{p + 2m}}}}} \right\}} {\zeta _n}\left( {1,2a + 1} \right)\\
&= \sum\limits_{i = 2}^{p + 2m - 1} {{{\left( { - 1} \right)}^{i - 1}}\zeta \left( {p + 2m + 1 - i,a + 1} \right)\left( {\sum\limits_{n = 1}^\infty  {\frac{{{\zeta _n}\left( {p,a + 1} \right)}}{{{{\left( {n + a} \right)}^i}}}} } \right)} \\
&\quad - \sum\limits_{i = 2}^{p - 1} {{{\left( { - 1} \right)}^{i - 1}}\zeta \left( {p + 1 - i,a + 1} \right)\left( {\sum\limits_{n = 1}^\infty  {\frac{{{\zeta _n}\left( {p + 2m,a + 1} \right)}}{{{{\left( {n + a} \right)}^i}}}} } \right)} \\
&\quad + \zeta \left( {p,a + 1} \right)\sum\limits_{n = 1}^\infty  {\frac{{{\zeta _n}\left( {1,a + 1} \right)}}{{{{\left( {n + a} \right)}^{p + 2m}}}}}  - \zeta \left( {p + 2m,a + 1} \right)\sum\limits_{n = 1}^\infty  {\frac{{{\zeta _n}\left( {1,a + 1} \right)}}{{{{\left( {n + a} \right)}^p}}}} \\
&\quad + \zeta \left( {p + 2m,a + 1} \right)\zeta \left( {p + 1,a + 1} \right) - \zeta \left( {p + 2m + 1,a + 1} \right)\zeta \left( {p,a + 1} \right)\\
&\quad + {\left( { - 1} \right)^{p - 1}}a\left( {\sum\limits_{n = 1}^\infty  {\frac{1}{{\left( {n + a} \right)\left( {n + 2a} \right)}}} } \right)\left( {\sum\limits_{n = 1}^\infty  {\frac{{{\zeta _n}\left( {p,a + 1} \right)}}{{{{\left( {n + a} \right)}^{p + 2m}}}}} } \right)\\
&\quad - {\left( { - 1} \right)^{p - 1}}a\left( {\sum\limits_{n = 1}^\infty  {\frac{1}{{\left( {n + a} \right)\left( {n + 2a} \right)}}} } \right)\sum\limits_{n = 1}^\infty  {\frac{{{\zeta _n}\left( {p + 2m,a + 1} \right)}}{{{{\left( {n + a} \right)}^p}}}},\tag{4.22}
\end{align*}
and
\begin{align*}
&{\left( { - 1} \right)^{p - 1}}\sum\limits_{n = 1}^\infty  {\left\{ {\frac{{{\zeta _n}\left( {p + 2m + 1,a + 1} \right)}}{{{{\left( {n + a} \right)}^p}}} + \frac{{{\zeta _n}\left( {p,a + 1} \right)}}{{{{\left( {n + a} \right)}^{p + 2m + 1}}}}} \right\}} {\zeta _n}\left( {1,2a+1} \right)\\
&= \sum\limits_{i = 2}^{p + 2m} {{{\left( { - 1} \right)}^{i - 1}}\zeta \left( {p + 2m + 2 - i,a + 1} \right)\left( {\sum\limits_{n = 1}^\infty  {\frac{{{\zeta _n}\left( {p,a + 1} \right)}}{{{{\left( {n + a} \right)}^i}}}} } \right)} \\
&\quad - \sum\limits_{i = 2}^{p - 1} {{{\left( { - 1} \right)}^{i - 1}}\zeta \left( {p + 1 - i,a + 1} \right)\left( {\sum\limits_{n = 1}^\infty  {\frac{{{\zeta _n}\left( {p + 2m + 1,a + 1} \right)}}{{{{\left( {n + a} \right)}^i}}}} } \right)} \\
&\quad + \zeta \left( {p,a + 1} \right)\sum\limits_{n = 1}^\infty  {\frac{{{\zeta _n}\left( {1,a + 1} \right)}}{{{{\left( {n + a} \right)}^{p + 2m + 1}}}}}  - \zeta \left( {p + 2m + 1,a + 1} \right)\sum\limits_{n = 1}^\infty  {\frac{{{\zeta _n}\left( {1,a + 1} \right)}}{{{{\left( {n + a} \right)}^p}}}} \\
&\quad + \zeta \left( {p + 2m + 1,a + 1} \right)\zeta \left( {p + 1,a + 1} \right) - \zeta \left( {p + 2m + 2,a + 1} \right)\zeta \left( {p,a + 1} \right)\\
&\quad + {\left( { - 1} \right)^p}a\left( {\sum\limits_{n = 1}^\infty  {\frac{1}{{\left( {n + a} \right)\left( {n + 2a} \right)}}} } \right)\left( {\sum\limits_{n = 1}^\infty  {\frac{{{\zeta _n}\left( {p,a + 1} \right)}}{{{{\left( {n + a} \right)}^{p + 2m + 1}}}}} } \right)\\
&\quad + {\left( { - 1} \right)^p}a\left( {\sum\limits_{n = 1}^\infty  {\frac{1}{{\left( {n + a} \right)\left( {n + 2a} \right)}}} } \right)\left( {\sum\limits_{n = 1}^\infty  {\frac{{{\zeta _n}\left( {p + 2m + 1,a + 1} \right)}}{{{{\left( {n + a} \right)}^p}}}} } \right).\tag{4.23}
\end{align*}
\end{cor}
Furthermore, setting $a=0$ in (4.22) and (4.23), we can give the following corollaries:
\begin{cor}
For $p\geq 2,m\geq 0$ and $p,m\in \N$, we have
\begin{align*}
&{\left( { - 1} \right)^{p - 1}}\sum\limits_{n = 1}^\infty  {\left\{ {\frac{{{H_n}{\zeta _n}\left( {p + 2m} \right)}}{{{n^p}}} - \frac{{{H_n}{\zeta _n}\left( p \right)}}{{{n^{p + 2m}}}}} \right\}}  \\
& =\sum\limits_{i = 2}^{p + 2m - 1} {{{\left( { - 1} \right)}^{i - 1}}\zeta \left( {p + 2m + 1 - i} \right)} \sum\limits_{n = 1}^\infty  {\frac{{{\zeta _n}\left( p \right)}}{{{n^i}}}}  - \sum\limits_{i = 2}^{p - 1} {{{\left( { - 1} \right)}^{i - 1}}\zeta \left( {p + 1 - i} \right)} \sum\limits_{n = 1}^\infty  {\frac{{{\zeta _n}\left( {p + 2m} \right)}}{{{n^i}}}}\\
&\quad + \zeta \left( p \right)\sum\limits_{n = 1}^\infty  {\frac{{{H_n}}}{{{n^{p + 2m}}}}}  - \zeta \left( {p + 2m} \right)\sum\limits_{n = 1}^\infty  {\frac{{{H_n}}}{{{n^p}}}}  + \zeta \left( {p + 1} \right)\zeta \left( {p + 2m} \right) - \zeta \left( p \right)\zeta \left( {p + 2m + 1} \right) .\tag{4.24}
\end{align*}
\end{cor}
\begin{cor}
Let $p \geq 2,m \ge 0$  be integers, we have
\begin{align*}
&{\left( { - 1} \right)^{p - 1}}\sum\limits_{n = 1}^\infty  {\left\{ {\frac{{{H_n}{\zeta _n}\left( {p + 2m + 1} \right)}}{{{n^p}}} + \frac{{{H_n}{\zeta _n}\left( p \right)}}{{{n^{p + 2m + 1}}}}} \right\}} \\
&=\zeta \left( p \right)\sum\limits_{n = 1}^\infty  {\frac{{{H_n}}}{{{n^{p + 2m + 1}}}}}  - \zeta \left( {p + 2m + 1} \right)\sum\limits_{n = 1}^\infty  {\frac{{{H_n}}}{{{n^p}}}}  + \zeta \left( {p + 1} \right)\zeta \left( {p + 2m + 1} \right) - \zeta \left( p \right)\zeta \left( {p + 2m + 2} \right)
\nonumber \\ &\quad + \sum\limits_{i = 2}^{p + 2m} {{{\left( { - 1} \right)}^{i - 1}}\zeta \left( {p + 2m + 2 - i} \right)} \sum\limits_{n = 1}^\infty  {\frac{{{\zeta _n}\left( p \right)}}{{{n^i}}}}  - \sum\limits_{i = 2}^{p - 1} {{{\left( { - 1} \right)}^{i - 1}}\zeta \left( {p + 1 - i} \right)} \sum\limits_{n = 1}^\infty  {\frac{{{\zeta _n}\left( {p + 2m + 1} \right)}}{{{n^i}}}} . \tag{4.25}
\end{align*}
\end{cor}
{\bf Acknowledgments.} The authors would like to thank the anonymous
referee for his/her helpful comments, which improve the presentation
of the paper.


{\small
\begin{thebibliography}{99}

\bibitem{BBG1994}
David H. Bailey, Jonathan M. Borwein and Roland Girgensohn. {\sl Experimental evaluation of Euler sums}.
Experimental Mathematics., 1994, {\bf 3}(1): 17-30.
\bibitem{BD2006}
Jonathan M. Borwein and David M. Bradley. {\sl Thirty-Two goldbach variations}. International Journal of Number Theory., 2006, {\bf 2}(1): 65-103.
\bibitem{BBD2008} D. Borwein, J. M. Borwein, and  D. M. Bradley. {\sl Parametric Euler sum identities}. Journal of Mathematical Analysis and Applications., 2008, {\bf 316}(1): 328-338.
\bibitem{BBG1995}
David Borwein, Jonathan M. Borwein and Roland Girgensohn. {\sl Explicit evaluation of
Euler sums}. Proc. Edinburgh Math., 1995, {\bf 38}: 277-294.
\bibitem{BBBL2001}
Jonathan M. Borwein, David M. Bradley, David J. Broadhurst, Petr. Lison¨§k.
{\sl Special values of multiple polylogarithms.} Trans. Amer. Math. Soc., 2001, {\bf 353}(3): 907-941.
\bibitem{BG1996}
J.M. Borwein, R. Girgensohn. {\sl Evaluation of triple Euler sums}, Electron. J. Combin., 1996: 2-7.
\bibitem{Boy2001}
K. Boyadzhiev. {\sl Evaluation of Euler-Zagier Sums, Internat}. J. Math. Sci., 2001, {\bf 27}(7): 407-412.
\bibitem{E2012}
Minking Eie, Chuan-Sheng Wei. {\sl Evaluations of some quadruple Euler sums of even weight}. Functions et Approximatio., 2012, {\bf 46}(1): 63-67.
\bibitem{FS1998}
Philippe Flajolet, Bruno Salvy. {\sl Euler sums and contour integral representations}. Experimental Mathematics., 1998, {\bf 7}(1): 15--35.
\bibitem{Fr2005}
Pedro freitas. {\sl Integrals of Polylogarithmic Functions, Recurrence Relations, and Associated Euler Sums}. Mathematics of Computation., 2005, {\bf 74}(251): 1425-1440.
\bibitem{MEH1992}
M. E. Hoffman. {\sl Multiple harmonic series}. Pacific J. Math., 1992, {\bf 152}: 275-290.
\bibitem{L1974}
 Comtet L. Advanced combinatorics, D Reidel Publishing Company, Boston,1974.
\bibitem{M2014}
I. Mez$\ddot{o}$. {\sl Nonlinear Euler sums}. Pacific J. Math., 2014, {\bf 272}: 201-226.
\bibitem{Ni1965}
N. Nielsen. {\sl Handbuch der Theorie der Gammafunction and Theorie des Integrallogarithmus and verwandter Transzendenten}., 1906, Reprinted together as Die Gammafunction, Chelsea, New York, 1965.
\bibitem {S2015}
A. Sofo. {\sl Quadratic alternating harmonic number sums}. J. Number Theory., 2015, {\bf 154}: 144-159.
\bibitem{S2011}
 A. Sofo, H.M. Srivastava. {\sl Identities for the harmonic numbers and binomial coefficients}. Ramanu-
jan J., 2011, {\bf 25}: 93-113.
\bibitem{S2010}
 A. Sofo. {\sl Harmonic sums and integral representations}. J. Appl. Anal., 2010, {\bf 16}: 265-277.
\bibitem{XC2016}
 Ce Xu, Jinfa Cheng. {\sl Some Results On Euler Sums}. Functions et Approximatio., 2016, {\bf 54}(1): 25-37.
\bibitem{X2016}
Ce Xu, Yuhuan Yan, Zhijuan Shi. {\sl Euler sums and integrals of polylogarithm functions}. J. Number Theory., 2016, {\bf 165}: 84-108.
\bibitem{CX2016}
Ce Xu, Mingyu Zhang, Weixia Zhu. {\sl Some evaluation of q-analogues of Euler sums}. Monatshefte F$\ddot{u}$r Mathematik ., 2016: 1-19.

\end{thebibliography}}
 \end{document}